\documentclass[11pt]{article}
\usepackage{latexsym,color,amsmath,amsthm,amssymb,amscd,amsfonts}

\newtheorem{lemma}{Lemma}[section]
\newtheorem{proposition}[lemma]{Proposition}

\newtheorem*{remark*}{Remark}

\def\R{{\mathbb R}}

\def\upddots{\mathinner{\mkern 1mu\raise 1pt \hbox{.}\mkern 2mu
\mkern 2mu \raise 4pt\hbox{.}\mkern 1mu \raise 7pt\vbox {\kern 7
pt\hbox{.}}} }

\newcommand{\spn}{{Sp_{2n}(\F)}}
\newcommand{\gsppn}{{GSp^{^+}_{2n}(\F)}}

\newcommand{\gspn}{{GSp_{2n}(\F)}}

\newcommand{\mspn}{{\overline{Sp_{2n}(\F)}}}

\newcommand{\tgspn}{T'_{2n}(\F)}

\newcommand{\bspn}{B_{2n}(\F)}
\newcommand{\bgspn}{B'_{2n}(\F)}
\newcommand{\bgsppn}{B^{^+}_{2n}(\F)}

\newcommand{\zspn}{{{N_{2n}}}(\F)}

\newcommand{\zgln}{N_{_{GL_n}}(\F)}

\newcommand{\mgspn}{\overline{GSp_{2n}(\F)}}

\newcommand{\mspr}{\overline{Sp_{2r}(\F)}}

\newcommand{\F}{F}

\newcommand{\Of}{\mathbb O_{\F}}

\newcommand{\spx}{\spn}
\newcommand{\vl}{v_{\lambda}}
\newcommand{\f}{_{_\F}}

\newcommand{\mspx}{\overline{Sp(X)}}
\newcommand{\mgspx}{\overline{GSp(X)}}
\newcommand{\gspx}{{\gspn}}
\newcommand{\mpt}{\overline{P_{\bold{t}}(\F)}}

\newcommand{\mgpt}{\overline{P'_{\bold{t}}(\F)}}

\newcommand{\mgppt}{\overline{P^{^+}_{\bold{t}}(\F)}}
\newcommand{\zt}{Z_{\bold{t}}(\F)}
\newcommand{\ztt}{Z_{\bold{t_2}}(\F)}
\newcommand{\mmt}{\overline{M_{\bold{t}}(\F)}}
\newcommand{\mmto}{\overline{M_{\bold{t_1}}(\F)}}
\newcommand{\mmtt}{\overline{M_{\bold{t_2}}(\F)}}
\newcommand{\ut}{U_{\bold{t}}(\F)}
\newcommand{\mgpspn}{\overline{GSp^{^+}_{2n}(\F)}}
\newcommand{\mgmt}{{\overline{M'_{\bold{t}}(\F)}}}
\newcommand{\mgmto}{{\overline{M'_{\bold{t_1}}(\F)}}}
\newcommand{\mgmtt}{{\overline{M'_{\bold{t_2}}(\F)}}}
\newcommand{\mgpmto}{{\overline{M^{^+}_{\bold{t_1}}(\F)}}}
\newcommand{\mgpmtt}{{\overline{M^{^+}_{\bold{t_2}}(\F)}}}

\newcommand{\gmto}{{M'_{\bold{t_1}}(\F)}}
\newcommand{\gmtt}{{M'_{\bold{t_2}}(\F)}}

\newcommand{\gpmtt}{{M^{^+}_{\bold{t_2}}(\F)}}
\newcommand{\mto}{{M_{\bold{t_1}}(\F)}}
\newcommand{\mtt}{{M_{\bold{t_2}}(\F)}}

\newcommand{\gpto}{{P'_{\bold{t_1}}(\F)}}

\newcommand{\gppto}{{P^{^+}_{\bold{t_1}}(\F)}}

\newcommand{\pto}{{P_{\bold{t_1}}(\F)}}

\newcommand{\mgptot}{{\overline{P'_{\bold{t_1}, \bold{t_2}}(\F)}}}
\newcommand{\mgpptot}{{\overline{P'^{^+}_{\bold{t_1}, \bold{t_2}}(\F)}}}
\newcommand{\mptot}{{\overline{P_{\bold{t_1}, \bold{t_2}}(\F)}}}
\newcommand{\gptot}{{P'_{\bold{t_1}, \bold{t_2}}(\F)}}
\newcommand{\gpptot}{{P'^{^+}_{\bold{t_1}, \bold{t_2}}(\F)}}
\newcommand{\ptot}{{P_{\bold{t_1}, \bold{t_2}}(\F)}}

\newcommand{\mt}{M_{\bold{t}}(\F)}
\newcommand{\gmt}{{M'_{\bold{t}}(\F)}}
\newcommand{\gpmt}{{M^{^+}_{\bold{t}}(\F)}}
\newcommand{\mgpmt}{{\overline{{M^{^+}_{\bold{t}}(\F)}}}}

\newcommand{\nt}{N_{\bold{t}}(\F)}

\newcommand{\pt}{{P_{\bold{t}}}(\F)}
\newcommand{\gpt}{{P'_{\bold{t}}}(\F)}
\newcommand{\gppt}{{P^{^+}_{\bold{t}}}(\F)}

\newcommand{\N}{\mathbb N}

\newcommand{\half}{\frac{1}{2}}

\newcommand{\ab} {|\!|}
\newcommand{\gl}{{GL_n(\F)}}
\newcommand{\gln}{{GL_n(\F)}}

\newcommand{\Q}{\mathbb Q}
\newcommand{\C}{\mathbb C}
\newcommand{\wc}{\widetilde{c}}

\def\>{\rangle}
\def\<{\langle}

\newtheorem{lem}{Lemma}[section]
\newtheorem{thm}{Theorem}[section]
\newtheorem{cor}{Corollary}[section]

\numberwithin{equation}{section}

\newcommand{\msl}{\overline{SL_2(\F)}}

\def\dotunion{
\def\dotunionD{\bigcup\kern-9pt\cdot\kern5pt}
\def\dotunionT{\bigcup\kern-7.5pt\cdot\kern3.5pt}
\mathop{\mathchoice{\dotunionD}{\dotunionT}{}{}}} 
\begin {document}

\author{Dani Szpruch}
\date{}
\title {Some results in the theory of genuine representations of the metaplectic double cover of $\gspn$ over p-adic fields.} \maketitle


\begin{abstract}
Let $\F$ be a p-adic field and let  $\overline{G(n)}$ and  $\overline{G_{_0}(n)}$ be the metaplectic double covers of the general symplectic group and symplectic group attached to a $2n$ dimensional symplectic space over $\F$. We show here that if $n$ is odd then all the genuine irreducible representations of $\overline{G(n)}$ are induced from a normal subgroup of finite index closely related to $\overline{G_{_0}(n)}$. Thus, we reduce, in this case, the theory of genuine admissible representations of $\overline{G(n)}$ to the better understood corresponding theory of $\overline{G_{_0}(n)}$. For odd $n$ we also prove the uniqueness of certain Whittaker functionals along with Rodier type of Heredity. Our results apply also to all parabolic subgroups of $\overline{G(n)}$ if $n$ is  odd and to some of the parabolic subgroups of $\overline{G(n)}$ if $n$ is even. We prove some irreducibility criteria for parabolic induction on $\overline{G(n)}$ for both even and odd $n$. As a corollary we show, among other results, that while for odd $n$, all genuine principal series representations of $\overline{G(n)}$ induced from unitary representations are irreducible, there exist reducibility points on the unitary axis if $n$ is even. We also list all the reducible genuine principal series representations of $\overline{G(2)}$ provided that the $\F$ is not 2-adic.
\end{abstract}
\section{Introduction}

In these notes we relate the theory of genuine admissible representations of the metaplectic double cover of $\gspn$, where $\F$ is a p-adic field and $n$ is odd,  to the corresponding theory of its derived group which is the metaplectic double cover of $\spn$. The prototype of our results are among the local results of Gelbart and Piatetski-Shapiro in \cite{GP80}, \cite{GP81} and \cite{GP}. In fact, for $n=1$ some of our main results were proven in these papers. Similar results were proven by Adams for the $m$ fold covers of $GL_m(\F)$ and $SL_m(\F)$, See \cite{A}. The authors of these papers used Clifford theory combined with explicit cocycle computations and we have followed their path. We also prove here the uniqueness of certain Whittaker functionals which is a generalization of the results of Gelbart, Howe and Piatetski-Shapiro, see \cite{GHP}. To explain our results we need some notation:

Let $G(n)=\gspn$ be the general symplectic group attached to a 2n dimensional symplectic space over $\F$, where $\F$ is a p-adic field. Let $G_{_0}(n)=\spn$ be the kernel of the similitude map and let $G^{^+}(n)$ be the subgroup of $G(n)$ which consists of elements whose similitude factor lies in ${\F^*}^2$. Denote by $\overline{G_{_0}(n)}$ the unique metaplectic double cover of $G_{_0}(n)$ and denote by $\overline{G(n)}$ the unique metaplectic double cover of $G(n)$ which contains $\overline{G_{_0}(n)}$. For a subset $A$ of $G$ we denote by $\overline{A}$ its primage in $\overline{G(n)}$. Let $Irr\bigl(\overline{G(n)}\bigr)$, $Irr\bigl(\overline{G^{^+}(n)}\bigr)$ and $Irr\bigl(\overline{G_{_0}(n)}\bigr)$ be the set of genuine smooth admissible irreducible representations of $\overline{G(n)}$, $\overline{G^{^+}(n)}$ and $\overline{G_{_0}(n)}$ respectively. For a group $H$ denote its center by $Z(H)$. It turns out that $\overline{G^{^+}}=\overline{G_{_0}}Z \bigl(\overline{G^{^+}(n)}\bigr)$ and that $\overline{G_{_0}} \cap Z \bigl(\overline{G^{^+}(n)}\bigr)$ is $\overline{\pm I_{2n}}$. Thus, given the representation theory of $\overline{G_{_0}}$, the representation theory of $\overline{G^{^+}(n)}$ is trivial. On the other hand our first main result is the following:

{\it{\bf Theorem A.}  Assume that $n$ is odd. The map
$$\pi \mapsto Ind_{\overline{G^{^+}(n)}}^{\overline{G(n)}}\pi$$ is a surjective map from $Irr \bigl(\overline{G^{^+}(n)} \bigr)$ to
$Irr \bigl(\overline{G(n)} \bigr)$.
Furthermore, any $\tau \in Irr \bigl(\overline{G(n)} \bigr)$ decomposes over $\overline{G^{^+}(n)}$ into a direct sum of $[{\F^*}:{\F^*}^2]$ pairwise non-isomorphic $\overline{G^{^+}(n)}$ modules. The inverse image of $\tau$ under the surjection mentioned above is the set of these summands.}

This result gives a natural correspondence between $Irr \bigl(\overline{G(n)} \bigr)$ and $Irr \bigl(\overline{G_{_0}} \bigr)$ which has no linear analog: Any $\tau \in Irr \bigl(\overline{G(n)} \bigr)$ decomposes over $\overline{G_{_0}(n)}$ into a direct sum of $[{\F^*}:{\F^*}^2]$ irreducible modules. We also show that unlike the linear case, multiplicity one does not hold. The importance of Theorem A lies in the fact that much is already known about the theory of genuine admissible representation of $\overline{G_{_0}(n)}$. An important feature of our result is that it continues to hold if $G(n)$ is replaced with any of its Levi subgroups. In fact, it holds for any Levi subgroup of $G(n)$ regardless of the parity of $n$, provided that the Levi subgroup in discussion is isomorphic to $$GL_{n_1}(\F) \times GL_{n_2}(\F) \ldots \times GL_{n_r}(\F) \times G(2n_{r+1})$$ where at least one of $n_1,\ldots, n_{r+1}$ is odd. We call these Levi subgroups, Levi subgroups of odd type and we call the other Levi subgroups, Levi subgroups of even type. Note that if $n$ is odd then all  Lvei subgroups of $G(n)$ are odd and that regardless of the parity of $n$, the group of diagonal matrices inside $G(n)$ is always of odd type. The upshot here is the following:

{\it{\bf Theorem B.}  Let $P$ be a parabolic subgroup of $G(n)$ and let $\pi$ be a genuine smooth admissible irreducible representation of $\overline{P}$. Define $P_{_0}=P \cap G_{_0}(n)$ and let $\pi_{_0}$ be any irreducible $\overline{P_{_0}}$ module which appears in $\pi$. Define $I(\pi)=Ind_{\overline{P}}^{\overline{G(n)}}\pi$ and  $I(\pi_{_0})=Ind_{\overline{P_{_0}}}^{\overline{G_{_0}(n)}}\pi_{_0}.$\\
1. If $n$ is odd then $I(\pi)$ is irreducible if and only if $I(\pi_{_0})$ is irreducible.\\
2. If $n$ is even, $P$ is of odd type and $\pi$ is supercuspidal then $I(\pi)$ is irreducible if and only if $I(\pi_{_0})$ is irreducible and $\pi_{_0}$ is not Weyl conjugate to any of its twists by elements of $\overline{P}$ which lies out side $\overline{G^{^+}(n)}$.}

{\it{\bf Theorem C.} Let $\pi$ be an element in $Irr \bigl(\overline{G(n)} \bigr)$. If either $n$ is odd or $\pi$ is induced from a parabolic subgroup of odd type then  its contragredient representation is isomorphic to a particular one dimensional non-genuine twist of $\pi$. For details see Theorem \ref{odd dual} and Proposition \ref{even dual}.}

Theorem B is still valid when $G$ is replaced with any parabolic subgroup $P'$ whose Levi part contains the Levi part of $P$ provided that $P$ of odd type. Our result implies that the number of reducibility points of $Ind_{\overline{P}}^{\overline{G(n)}}\pi$ is no less then the number of reducibility points of $Ind_{\overline{P_{_0}}}^{\overline{G_{_0}(n)}}\pi_{_0}$. It is well known that for linear groups, the opposite holds.
Using some results from \cite{Sz12} and Theorem B we prove some irreducibility theorems on genuine parabolic induction on $\overline{G(n)}$. In particular we prove the following:

{\it {\bf Theorem D.} If $n$ is odd then all genuine principle series representations of $\overline{G(n)}$ induced from unitary representations are irreducible. If $n$ is even then there exist reducible genuine principle series representations of $\overline{G(n)}$ induced from unitary representations. Among these reducible reorientations are also unramified representations.}

Note that, regardless of the parity of $n$, all genuine principle series representations of $\overline{G_{_0}(n)}$ induced from unitary characters are irreducible. Using the results of Zorn, \cite{Zo09}, we also list all the reducible principle series representations of $\overline{G(2)}$ over p-adic fields of odd residual characteristic.

Assume now that $P_{_0}$ is the Siegel Parabolic subgroup of $G_{_0}(n)$. Let $\pi_{_0}$ be an irreducible smooth admissible  generic representation of $\overline{P_{_0}}$. The Plancherel measure, $\mu(\tau,s)$, attached to
$Ind^{\overline{G_{_0}(n)}}_{\overline{P_{_0}}} \pi_{_0}$
was computed in Theorem 4.3 of \cite{Sz12}. It was shown there that if $\chi$ is a quadratic character of $\F^*$ then $\mu(\pi_{_0},s)$ and $\mu(\chi \! \circ \! \det \otimes \pi_{_0},s)$ have the same analytic properties. This fact, which has no linear analog, is explained here by the observation that these two induced representations are conjugate by an element of $\overline{G(n)}$.

We now state our result on Whittaker functionals: If $n$ is odd then $$[Z\bigl(\overline{G^{^+}(n)}\bigl):Z \bigl(\overline{G(n)}\bigr)]=[{\F^*}:{\F^*}^2].$$
For $\pi$, a representation of $\overline{G(n)}$ where $n$ is odd, with a central character $\chi_\pi$, let $\Omega_\pi$ denote the finite set of characters of $Z \bigl(\overline{G^{^+}(n)}\bigr)$ which extends $\chi_\pi$. Let $N(n)$ ba a maximal unipotent radical of $G(n)$ and let $\psi$ be a non-degenerate character of $N(n)$. Let $W_\psi$ be the space of $\psi$-Whittaker functionals on $\pi$. For $\omega \in \Omega_\pi$, let $W_{\psi \times \omega}$ be the subspace of $W_\psi$ which consists of eigen functionals corresponding to $\omega$.

{\it{\bf Theorem E.} Let $\pi$ be smooth admissible generic representation of $\overline{G(n)}$ where $n$ is odd. Assume that $\pi$ is either irreducible or parabolically induced from a smooth admissible irreducible representation. Then,  for any $\omega \in \Omega_\pi$, $\dim(W_{\psi \times \omega}) \leq 1$ and $\dim(W_{\psi \times \omega}) = 1$ for at least one $\omega \in \Omega_\pi$.} Furthermore, if $n$ is even and $\pi$ is  parabolically induced from a smooth admissible irreducible representation of an odd parabolic subgroup then one may generalize this uniqueness and exitance. For details see Theorem \ref{rodier gen}.

This result is proven by combining Theorem A and the Uniqueness of $\psi$-Whittaker functional for $\overline{G_{_0}(n)}$ proven in \cite{Sz}. It should be noted that this is not the argument used in \cite{GHP} for $n=1$. From Theorem E it follows that $\dim(W_{\pi} ) \leq [\F^*:{\F^*}^2]$. In some case sharper results are obtained. We now explain some of the expected applications of the results given in this paper:

1. Theorem E enables the definition of an analog to Shahidi local coefficients, \cite{Sha 1}, for a parabolic induction on $\overline{G(n)}$, provided that the inducing subgroup is of odd type: $\omega\times \psi$ and its generalization mentioned in Theorem E replace the role of the usual $\psi$. The computation of these coefficients is reduced at once to $\overline{G_{_0}(n)}$. One can use these local coefficients to define local $\gamma,\, L$ and $\epsilon$ factors in a way which is similar to Shahidi's definitions in \cite{Sha90}. It can be proven that these factors satisfies cerian properties that characterize them uniquely, see Theorem 3.5 of \cite{Sha90} and Theorem 4 of \cite{L}. For instance, the multiplicativity of the gamma factor will follow from the multiplicativity of the local coefficients along with the fact that a parabolic subgroup of odd Levi subgroup is always of odd type.

2. In the minimal parabolic case, our results are in accordance with the construction of principal series for covering groups, see \cite{Mc}. It turns out that if $T(n)$ is the diagonal subgroup of $G(n)$ then $\overline{T \cap G^{^+}(n)}$ is a maximal abelian subgroup of $\overline{T(n)}$. Note that it is not an analog to Kazhdan Patterson standard maximal Abelian subgroup, see Section 1 of \cite{KP}, which appears also in \cite{CO} and \cite{Mc2}. Its important property is that it contains the commutative group $\overline{T \cap G^{^0}(n)}$. In the unramified case, one may use Theorem E construct a basis for the space of spherical Whittaker functions which consists of symmetric functions and whose functional equation is diagonal. Given the global functional equation satisfied by Eisenstein series, see Theorem IV.1.10 of \cite{MW95}, this fact is meaningful. An explicit formulas here may be given using the results of Bump, Friedberg and Hoffstein in \cite{BFH}.

3. Let $SO_{2n+1}^{\pm}(\F)$ be the two special orthogonal groups corresponding to the two $2n+1$ orthogonal spaces over $\F$
with discriminant equals 1. In \cite{GS}, Gan and Savin established a bijection between $Irr \bigl(\overline{G_{_0}(n)} \bigr)$ and the irreducible smooth admissible representations of $SO_{2n+1}^{\pm}(\F)$ closely related to the theta correspondence, and used it to prove the local Langlands conjecture for $\overline{G_{_0}(n)}$ given the local Langlands conjecture for $SO_{2n+1}^{\pm}(\F)$. Using the frame work of Roberts, \cite{Ro}, this correspondence extends easily to a correspondence between the irreducible smooth admissible representations of $GSO_{2n+1}^{\pm}(\F)$ and $Irr \bigl(\overline{G^{^+}(n)}\bigr)$. Assume now that $n$ is odd. Using Theorem A, it is possible to define a natural correspondence between irreducible smooth admissible representations of $GSO_{2n+1}^{\pm}(\F)$ and $Irr \bigl(\overline{G(n)} \bigr)$. This correspondence will preserve local factors and Plancherel measures, at least for generic representations. This correspondence can then be used to prove the local Langlands conjecture for $\overline{G(n)}$ given the local Langlands Conjecture for  $GSO_{2n+1}^{\pm}(\F)$.

We shall address each of these topics in a future publication.

This paper is organized as follows: In section \ref{preliminaries} we introduce some notation and recall some facts about the metaplectic groups and Rao's cocycle. In Section \ref{structural facts} we carry out the crucial cocycle computation and prove some structural facts about $\overline{G(n)}$  and its subgroups. We use these facts in Section \ref{rrt} where we apply Clifford theory to prove Theorems A-C and other related results. The irreducibility theorems are proven in Section \ref{irr thm} using the results of Section \ref{rrt}. In Section  \ref{whi} we prove Theorem E and other results using the results of Section \ref{rrt} and uniqueness of Whittaker model for $\overline{G_{_0}(n)}$.

I would like to thank Freydoon Shahidi, Gordan Savin, Jeffery Adams, Nadya Gurevich and  Sandeep Varma for useful communications on the subject matter.

\section{preliminaries} \label{preliminaries}
\subsection{General notations} \label{genral not}

Through these paper $\F$ will denote a p-adic field, i.e, a finite extension of $\Q _p$. If $G$ is a group we denote its center by $Z(G)$. For  $g,h \in G$ we define $$h^g=g^{-1}hg.$$ Let $\pi$ be a representation of $G$. If $\pi$ has a central character, we denote it by $\chi_\pi$. If $H$ is a normal subgroup of $G$, $\tau$ is a representation of $H$ and $g$ is an element of $G$, we denote by $\tau^g$ the $H$ module defined by $$h \mapsto \tau(h^g).$$
If $G$ is a p-adic group and $\pi$ is a smooth representation of $G$, we denote by $\widehat{\pi}$ the contragredient representation.

\subsection{Hilbert symbol and Weil index} \label{hw}

Let $(\cdot,\cdot)\f$ be the quadratic Hilbert symbol of $\F$. It is a non-degenerate symmetric bilinear form on $\F^* /{\F^*}^2$. Recall that for all $a,\, b \in \F^*$
\begin{equation} \label{Hilbert properties}
(a,b)\f=(a,-ab)\f.
\end{equation}
For $a \in \F^*$ we define $\eta_a$ to be the quadratic character of $\F^*$ attached to $a$, that is, $$\eta_a(b)=(a,b)\f.$$
Recall that $\eta_a=\eta_{ab^2}$. Let $\psi$ be a non-trivial additive character of $\F$. For $a\in \F^*$ define  $$\psi_a(x)=\psi(ax).$$ It is also a non-trivial additive character of $\F$. For $a\in \F^*$ let
$\gamma_\psi(a) \in \C^{1}$ be the normalized Weil factor associated with
the character of second degree of $\F$ given by $x \mapsto
\psi_a(x^2)$ (see Theorem 2 of Section 14 of \cite{Weil}). It is
known that $\gamma_{\psi}$ is a forth root of unity and that $\gamma_{\psi}\bigl({\F^*}^2\bigr)=1$. Also,
\begin{equation} \label{gammaprop}
\gamma_{\psi}(ab)=\gamma_{\psi}(a)\gamma_{\psi}(b)(a,b)_{\F}, \, \,  \, \,   \gamma_{\psi_b}=\eta_b \cdot \gamma_\psi  ,
\, \,  \, \,  \gamma_{\psi}^{-1}=\gamma_{\psi^{-1}}.
\end{equation}
\subsection{Linear groups} \label{The symplectic group}
Let $\gspn$ be the general symplectic group attached to a 2n dimensional symplectic space over $\F$. We shall realize $\gspn$ as the group $$\{g \in GL_{2n}(\F) \mid gJ_{2n}g^t=\lambda(g)J_{2n} \},$$
where
$$J_{2n}=\begin{pmatrix} _{0} & _{I_{n}}\\_{-I_{n} } & _{0}
\end{pmatrix}$$ and $\lambda(g)\in \F^*$ is the similitude factor of $g$. The similitude map $g\mapsto \lambda(g)$ is a rational character on $\gspn$. The kernel of the similitude map is
is the symplectic group, $\spn$. $\F^*$ is embedded in $\gspx$ via $$\lambda
\mapsto i(\lambda)=\begin{pmatrix} _{I_n} & _{0}\\_{0} & _{\lambda
I_n}\end{pmatrix}.$$ Using this embedding we define an action of
$\F^*$ on $\spx$: $$(g,\lambda) \mapsto
g^{i(\lambda)}.$$ Let $\F^* \ltimes
\spx$ be  the semi-direct product corresponding to this action.
For $g \in \gspx$ define $$g_{_1}=i\bigl(\lambda^{-1}(g)\bigr)g \in \! \spx.$$
 The map $$g \mapsto \bigl(\lambda(g),g_{_1} \bigr)$$ is an isomorphism between
$\gspx$ and $\F^* \ltimes \spx.$ In particular $GSp_0(\F) \simeq \F^*$.

We define $\gsppn$ to be the subgroup of $\gspn$ which consists of elements
whose similitude factor lies in ${\F^*}^2$. $\gsppn$ is a normal subgroup of $\gspn$ which contains $\spn$ . Clearly $$[\gspn:\gsppn]=[{\F^*}:{\F^*}^2] < \infty.$$
For $0 \leq r \leq n$ define $i_{r,n}$ to
be an embedding of $GSp_{2r}(\F)$ in $GSp_{2n}(\F)$ by $$
g=\begin{pmatrix} _{a} & _{b}\\_{c} & _{d}
\end{pmatrix}\mapsto \begin{pmatrix} _{I_{n-r}} & _{ } & _{ }
\\ _{ } & _{a} & _{ } & _{b}\\ _{ } & _{ }  & _{\lambda(g)I_{n-r}} & _{ } \\
_{ } & _{c}  & _{ } & _{d}
\end{pmatrix},$$ where $a,b,c,d \in Mat_{r \times r}(\F)$. The restrictions of $i_{r,n}$ and to $GSp^{^+}_{2r}(\F)$ and to $Sp_{2r}(\F)$ are embeddings of $GSp^{^+}_{2r}(\F)$ and $Sp_{2r}(\F)$ inside $\gsppn$  and $\spn$ respectively.

Let $\zgln$ be the group of upper triangular unipotent matrices in
$\gl$ and let $\zspn$ be the following maximal
unipotent subgroup of $\spn$:
$$\Bigl\{\begin{pmatrix} _{z} & _{b}\\_{0} & _{\widetilde{z}}
\end{pmatrix}  \mid z\in \zgln, b \in Mat_{n \times n}
(\F), \, b^t=z^{-1}b z^t \Bigr \},$$ where for $a \in GL_n$ we
define $\widetilde{a}={^t \!{a} \!^ {-1}}$. Let $\tgspn$ be the subgroup of diagonal matrices
inside $\gspn$. Denote $$\bgspn=\tgspn \ltimes \zspn.$$
It is a Borel subgroup of $\gspn$. A parabolic subgroup of $\gspn$ is called standard if it contains $\bgspn$. A standard Levi (unipotent) subgroup is a Levi (unipotent) part of a standard parabolic subgroup. Let $n_1,n_2,\ldots,n_r,n_{r+1}$ be $r+1$ nonnegative integers whose sum
is $n$. Put $$\bold{t}=(n_1,n_2,\ldots,n_r;n_{r+1}).$$ Let
$\gmt$ be the standard Levi subgroup of $\gspn$
which consists of elements of the form
$$[g_{_1},g_2,\ldots,g_r;h]=
diag(g_{_1},g_2,\ldots,g_r,I_k,\lambda(h)\widetilde{g_{_1}},\lambda(h)\widetilde{g_2},\ldots,\lambda(h)\widetilde{g_r},I_k)i_{k,n}(h),$$
where $g_i \in GL_{n_i}(\F), h \in GSp_{2n_{r+1}}(\F)$. Define
\begin{eqnarray} \nonumber \gpmt \! \! \! &=& \! \! \! \gmt\cap \gsppn, \\ \nonumber \mt\! \! \! &=& \! \! \!\gmt \cap \spn.\end{eqnarray}
Note that
\begin{eqnarray} \nonumber \gmt \! \! \! &\simeq& \! \! \! GL_{n_1}(\F) \times GL_{n_2}(\F) \ldots \times GL_{n_r}(\F) \times GSp_{2n_{r+1}}(\F),
\\ \nonumber \mt \! \! \! &\simeq& \! \! \! GL_{n_1}(\F) \times GL_{n_2}(\F) \ldots \times GL_{n_r}(\F) \times Sp_{2n_{r+1}}(\F),
\\ \nonumber \gpmt  \! \! \! &\simeq& \! \! \! GL_{n_1}(\F) \times GL_{n_2}(\F) \ldots \times GL_{n_r}(\F) \times GSp_{2n_{r+1}}^{^+}(\F).\end{eqnarray}
Define $$\gpt= \gmt \ltimes \ut$$ to be the standard parabolic subgroup of $\gspn$ whose Levi part is $\gmt$. Define now
\begin{eqnarray} \nonumber \pt  \! \! \! &=& \! \! \! \mt \ltimes \ut,
\\  \nonumber \gppt \! \! \! &=& \! \! \! \gpmt \ltimes \ut.\end{eqnarray}
$\pt$ and $\gppt$ are parabolic subgroups of $\spn$ and $\gspn$ respectively. We continue to call the parabolic and Levi subgroups obtained that way standard. To be clear we note here that by a parabolic subgroup we do not necessarily mean a proper parabolic subgroup. In fact, many of the results given in these notes for parabolic subgroups and their metaplectic double covers are motivated primarily by the $$\gmt=\gpt=\gspn$$ and $$ \! \! \! \! \! \mt=\pt=\spn$$ cases.

A particular role is reserved for $P_{(n;0)}(\F)$, the
Siegel parabolic subgroup of $\spn$:
\begin{equation}\label{parabolic} P_{(n;0)}(\F)= \Bigl\{
\begin{pmatrix} {a} & {b}\\{0} & {\widetilde{a}}
\end{pmatrix} \mid a \in \gl, b \in Mat_{n \times n} (\F),
b^{t}=a^{-1}ba^t \Bigr\}.\end{equation}
The Bruhat decomposition of $\spn$ with respect to $P_{(n;0)}(\F)$ is
$$\spn=\bigcup_ {j=0}^n \Omega_j (\F),$$
where $\Omega_j (\F)$, the $j^{th}$ Bruhat cell, is
$$\Omega_j(\F)=\bigl\{ \begin {pmatrix} _{\alpha} & _{\beta}\\_{\gamma } &
_{\delta}
\end{pmatrix} \in \spn \mid \alpha,\beta,\gamma, \delta \in Mat_{n \times n}(\F),\ rank(\gamma)=j \bigr\}.$$
Note that $\Omega_0 (\F)=P_{(n;0)}(\F)$, $\Omega_j^{i(\lambda)}(\F)=\Omega_j (\F)$ for all $\lambda \in \F^*$  and that $\Omega_j^{-1} (\F)=\Omega_j (\F)$. For $g=\begin{pmatrix} {a} & {b}\\{0} & {\widetilde{a}}
\end{pmatrix} \in \Omega_0(\F)$ define \begin{equation} \label{x def} x(g)=\det(a) \in \F^* / {\F^*}^2.\end{equation}
In Lemma 5.1 of \cite{R}, Rao extended the map  $g \mapsto x(g)$ to a certain map from $\spn$ to $\F^* / {\F^*}^2$
such that if $p_1, \, p_2 \in \Omega_0(\F)$ then
\begin{equation} x(p_1gp_2)=x(p_1)x(g)x(p_2)\end{equation}
for all $g \in \spn$. This map plays an important roll in the definition of Rao`s cocycle, see Section  \ref{des rao}. In page 362 of \cite{R} and (2-11) of page 457 of \cite{Sz} respectively, it is shown that for $g \in \Omega_j(\F)$ we have.
\begin{eqnarray} x(g^{-1}) &=& (-1)^j x(g),
\\ \label{x lambda} x(g^{i(\lambda)}) &=& \lambda^jx(g).\end{eqnarray}

\subsection{Rao's cocycle} \label{des rao}
Let $\mspn$ be the unique non-trivial two-fold cover of $\spn$. We shall realize $\mspn$ as the set
$\spx \times \{ \pm 1 \}$ equipped with the multiplication law
\begin{equation} \label{metaplectic structure}(g_{1},\epsilon_{1})(g_{2},\epsilon_{2})=
\bigl( g_{1}g_{2},\epsilon_{1}\epsilon_{2}c(g_{1},g_{2})\bigr)
.\end{equation} Here $$c:\spn \times \spn \rightarrow \{\pm 1\}$$ is Rao's cocycle, see \cite{R}.
For any subset $H$ of $\spn$ we denote by $\overline{H}$ its inverse image in $\mspn$ and we continue to denote by $H$ the subset $(H,1)$ of $\mspn$. We shall later use similar notations in the context of $\mgspn$. We recall the following explicit properties of Rao's cocycle. For $g\in \Omega_j(\F), \, p\in \Omega_0(\F)$ we have

\begin{eqnarray} \label{coc inv} c(g,g^{-1})&=&\bigl( x(g),(-1)^jx(g)\bigr)\f
(-1,-1)\f^{\frac{j(j-1)}{2}}
\\ \label{coc with p} c(p,g)&=&c(g,p)=\bigl(
x(p),x(g) \bigr)\f. \end{eqnarray}
See Corollaries 5.4 and 5.5 in \cite{R} (there
is a small mistake there which was corrected by
Adams in Theorem 3.1 of \cite{Kud}). A representation $\pi$ of $\mspn$ is called genuine if $$\pi(I_{2n},1)=-Id.$$ Similar definition holds for all the covering groups discussed in this paper.

\begin{lem}\label{mspn prop} The following hold:\\
1.  Two elements in $\mspn$  commute if and only if their projections to $\spn$ commute.\\
2. $\mspn$ splits over any standard unipotent subgroup of $\spn$ via the trivial section. Furthermore
$$\mpt= \overline{\mt} \ltimes \ut,$$
3. $(h,\epsilon) \mapsto (i_{r,n}(h),\epsilon)$ defines an embedding of $\mspr$ inside $\mspn$.\\
4. Fix $\bold{t}={(n_1,n_2,\ldots,n_r;n_{r+1})}$. For $1 \leq i \leq r$, let  $\sigma_i$ be a smooth irreducible representation of $GL_{n_i}(\F)$ and let $\overline{\sigma}$ be a representation of $\overline{Sp_{2n_{r+1}}(\F)}$. Fix $\psi$, a non-trivial additive character of $\F$, and define:

$$\tau\bigl([g_1,g_2,\ldots,g_{n_r};h],\epsilon \bigr) \mapsto \gamma_\psi\bigl(\prod_{i=1}^{r}\det(g_i)\bigr)\otimes \bigl(\otimes_{i=1}^r
{\sigma_i}(g_i)\bigr)\otimes \overline{\sigma}(h,\epsilon),$$
$\tau$ is an irreducible smooth admissible genuine representation of $\mmt$, and all  smooth admissible genuine irreducible representations of $\mmt$ are obtained that way.
 \end{lem}
\begin{proof} These are all well known: 1 follows from page 39 of \cite{MVW}. 2 follows from \eqref{coc inv} and \eqref{coc with p}. 3 follows from the inductive property of Rao's cocycle; See corollary 5.6 of \cite{R}. For 4 see Section 2.3 of \cite{GS} for example.
\end{proof}

Remark: It was proven in \cite{Sz phd} that if $\F$ is a p-adic field of odd residual characteristic then

$$ \mmt \simeq GL_{n_1}(\F) \times GL_{n_2}(\F) \ldots \times GL_{n_r}(\F)
\times \overline{Sp_{2n_{r+1}}(\F)} .$$
We shall make no direct use of this fact.

\subsection{Extension of Rao's cocycle to $\gspn$.} \label{ext rao}
In this section we shall recall the construction of  $$\mgspn \simeq \F^* \ltimes \mspn,$$ the unique two-fold cover of $\gspn$ which contains $\mspn$. By page 36 of \cite{MVW}, it is known that we can lift uniquely the outer conjugation $g \mapsto g^{i(\lambda)}$
of $\spn$ to $\mspn$, i.e, we can define
a map $v_\lambda:\spx \rightarrow \{ \pm 1 \}$ such that
$$(g,\epsilon) \mapsto
(g,\epsilon)^{\lambda}=\bigl(g^{\lambda},\epsilon
v_{\lambda}(g)\bigr)$$ is an automorphism of $\mspx$. In Section 2B of \cite{Sz} it was shown that
for $g \in \Omega_j(\F), \, \lambda \in \F^*,$ we have

\begin{equation} \label{the formula}\vl(g)=\bigl(x(g),\lambda^{j+1}\bigr)\f\bigl(\lambda,\lambda\bigr)\f^{\frac{j(j-1)}{2}}.\end{equation}
Furthermore, the map
$$\bigl(\lambda,(g,\epsilon)\bigr) \mapsto
(g,\epsilon)^{i(\lambda)}$$ defines an action of $\F^*$ on $\mspx$.
As explained in  Section 2B of \cite{Sz}, this computation enables the extension of $c(\cdot,\cdot)$ to a 2-cocycle
$\widetilde{c}(\cdot,\cdot)$ on $\gspx$ which takes values in
$\{\pm 1\}$: We define the group
$\F^* \ltimes \mspx$ using the multiplication formula
$$\bigl(a,(g,\epsilon_1)\bigr) \bigl
(b,(h,\epsilon_2)\bigr)=\bigl(ab,(g,\epsilon_1)^{i(b)}(h,\epsilon_2)\bigr).$$
We now define a bijection from the set $\mgspn=\gspx \times \{ \pm 1 \}$
to the set $\F^* \times \mspx$ by the formula
$$\overline{\iota}(g,\epsilon)=\bigl(\lambda(g),(g_{_1},\epsilon)\bigr),$$
whose inverse is given by
$$\overline{\iota}^{-1}\bigl(\lambda,(h,\epsilon)\bigr)=\bigl(i(\lambda)h,\epsilon
\bigr ).$$ We use $\overline{\iota}$ to define a group structure
on $\gspx \times \{ \pm 1 \}$. A straightforward computation shows
that the multiplication in $\mgspx$ is given by
$$(g,\epsilon_1)(h,\epsilon_2)=\bigl(gh,v_{\lambda(h)}(g_{_1})c \bigl(g_{_1}^{i(\lambda(h))},h_1\bigr)\epsilon_1\epsilon_2 \bigr).$$
Thus,
\begin{equation} \label{extension}
\widetilde{c}(g,h)=v_{\lambda(h)}\bigl(g_{_1} \bigr )c \bigl
(g_{_1}^{\lambda(h)},h_{_1} \bigr).
\end{equation}
We remark here that Kubota, \cite{Kub},
used a similar construction to extend a non-trivial double cover
of $SL_2(\F)$ to a non-trivial double cover of $GL_2(\F)$. For
$n=1$ our construction agrees with Kubota's. Note that by restricting $\widetilde{c}$ we obtain a double cover of $\gsppn$. Clearly
$$[\mgmt:\mgpmt]=[\gmt:\gpmt]=[\F^*:{\F^*}^2].$$
As a set of representatives of $\mgmt / \mgpmt$ we may take $$\{\bigl(i(\lambda),1\bigr)\mid \lambda \in \F^* /{\F^*}^2 \}.$$
We note that from the third assertion of Lemma \ref{mspn prop} and from \eqref{extension} it follows that  $$(h,\epsilon) \mapsto (i_{r,n}(h),\epsilon)$$ defines embeddings of $\overline{GSp^{^+}_{2r}(\F)}$ and $\overline{GSp_{2r}(\F)}$ inside $\mgpspn$  and $\mgspn$ respectively.

Notation:  For $\sigma=(g,\epsilon)$, a pullback in $\mgspn$ of an element $g\in \gspn$ we define $\lambda(\sigma)=\lambda(g)$.

We note here that fixing $(g,\epsilon) \in \mgspn$, the inner automorphism of $\mgspn$, $\sigma \mapsto \sigma^{(g,\epsilon)}$ does not depend on $\epsilon$. Thus, there is no ambiguity in the notation $\sigma \mapsto \sigma^g$ where $\sigma \in \mgspn,$ and $g \in \gspn$.

\section{Some structural facts about $\mgspn$ and its subgroups} \label{structural facts}

\begin{lem} \label{inver cw}Fix $g \in \gspn$. Assume that $g_{_1} \in \Omega_j(\F)$. Define $\lambda=\lambda(g)$. Then, the following generalization of \eqref{coc inv} holds:
$$\wc(g,g^{-1})=\bigl(-\lambda,-1\bigr)\f^{\frac{j(j-1)}{2}}\bigl(x(g_{_1}),(-\lambda)^{j+1}\bigr)\f$$

\end{lem}

\begin{proof} First note that since $g^{-1}=i(\lambda^{-1})\bigl((g_{_1})^{-1}\bigr)^{i(\lambda^{-1})}$ we have

\begin{equation} \label{pg pg-1} (g^{-1})_{_1}=\bigl((g_{_1})^{-1}\bigr)^{i(\lambda^{-1})}=\bigl((g_{_1})^{i(\lambda^{-1})}\bigr)^{-1} .\end{equation}
This implies that
\begin{eqnarray} \nonumber \wc(g,g^{-1})&=&v_{_{\lambda^{-1}}}\bigl(g_{_1} \bigr)c\bigl((g_{_1})^{i(\lambda^{-1})},(g^{-1})_{_1}\bigr) \\ \nonumber
&=&\bigl(x\bigl(g_{_1}\bigr),\lambda^{j+1}\bigr)\f\bigl(\lambda,\lambda\bigl)\f^{\frac{j(j-1)}{2}}
c\bigl((g_{_1})^{i(\lambda^{-1})},\bigl((g_{_1})^{i(\lambda^{-1})}\bigr)^{-1}\bigr).\end{eqnarray}
Hence, by \eqref{coc inv} we conclude that
$$\wc(g,g^{-1})= \bigl(x(g_{_1}),\lambda^{j+1}\bigr)\f\bigl(\lambda,\lambda\bigl)\f^{\frac{j(j-1)}{2}}
\bigl(x((g_{_1})^{i(\lambda^{-1})}),(-1)^jx((g_{_1})^{i(\lambda^{-1})})\bigr)\f
\bigl(-1,-1\bigr)\f^{\frac{j(j-1)}{2}}.$$
By \eqref{x lambda} and \eqref{Hilbert properties} we now get
\begin{eqnarray} \nonumber \wc(g,g^{-1}) &=& \bigl(-\lambda,-1 \bigr)\f^{\frac{j(j-1)}{2}}\bigl(x\bigl(g_{_1}\bigr),\lambda^{j+1}\bigr)\f
\bigl(x\bigl(g_{_1}\bigr)(\lambda)^j,x\bigl(g_{_1}\bigr)(-\lambda)^j\bigr)\f \\ \nonumber
&=&\bigl(-\lambda,-1 \bigr)\f^{\frac{j(j-1)}{2}}\bigl(x\bigl(g_{_1}\bigr),\lambda^{j+1}\bigr)\f
\bigl(x\bigl(g_{_1}\bigr),x\bigl(g_{_1}\bigr)(-1)^j\bigr)\f.\end{eqnarray}
Using \eqref{Hilbert properties} again the Lemma is now proven.

\end{proof}
\begin{lem} \label{ugly fact}Assume that $g \in \gspn$ is such that $g_{_1} \in \Omega_0(\F)$. Denote $\lambda=\lambda(g)$. For $h \in \gspn$ such that $h_1 \in \Omega_j(\F)$ we have
\begin{equation} \label{the fact}(g,\epsilon)(h,\epsilon')(g,\epsilon)^{-1}=\bigl(g h g^{-1},\epsilon' d(g,h) \bigr),\end{equation}
where \begin{equation} d(g,h)=\bigl(x(g_{_1}),\lambda(h)\bigr)\f\bigl(x(h_{_1}),\lambda^{j+1}\bigr)\f \bigl(\lambda,\lambda \bigr)\f^{\frac{j(j-1)}{2}}. \end{equation}

\end{lem}
\begin{proof} We have to show that \begin{equation} \label{d formula}\wc(g,h)\wc(g h,g^{-1})\wc(g,g^{-1})=d(g,h).\end{equation}
First, by Lemma \ref{inver cw} and by \eqref{Hilbert properties} we have
\begin{equation} \label{d formula 1}\wc(g,g^{-1})=\bigl(x(g_{_1}),-\lambda\bigr)\f.\end{equation}
Next, from \eqref{extension}, \eqref{x lambda} and \eqref{coc with p} it follows that
\begin{equation} \label{d formula 2}\wc(g,h)=\bigl(x(g_{_1}),\lambda(h)\bigr)\f\bigl(x(g_{_1}),x(h_{_1})\bigr)\f. \end{equation}
By \eqref{d formula}, \eqref{d formula 1} and \eqref{d formula 2}, the proof is done once we show that
\begin{equation} \label{d formula 3}\wc(gh,g^{-1})=\bigl(x(g_{_1})x(h_{_1}),\lambda^{j+1}\bigr)\f \bigl(\lambda,\lambda\bigr)\f^{\frac{j(j-1)}{2}}
\bigl(\lambda^j x(g_{_1})x(h_{_1}),x(g_{_1})\bigr)\f.\end{equation}
Note that $gh= i \bigl(\lambda \lambda(h)\bigr)(g_{_1})^{i(\lambda(h))}h_{_1}$. Therefore, $$(gh)_{_1}=(g_{_1})^{i(\lambda(h))}(h_{_1})\in \Omega_j(\F), \, \, x\bigl((gh)_{_1}\bigr)=x(g_{_1})x(h_{_1}\bigr).$$ Also, since $$\bigl ((gh)_{_1}\bigr)^{i(\lambda^{-1})}= (g_{_1})^{i(\lambda(h) \lambda^{-1})}(h_{_1})^{i(\lambda^{-1})},$$ we conclude that $\bigl((gh)_{_1}\bigr)^{i(\lambda^{-1})} \in \Omega_j(\F)$ and that $$x\bigl(((gh)_{_1})^{i(\lambda^{-1})}\bigr)=\lambda^j x(g_{_1})x(h_{_1}).$$
Thus,
\begin{eqnarray} \nonumber \wc(gh,g^{-1}) &=& v_{_\lambda}\bigl((gh)_{_1}\bigr)c\bigl(((gh)_{_1})^{i(\lambda^{-1})},(g^{-1})_{_1}\bigr) \\ \nonumber
&=&\bigl(x(g_{_1})x(h_{_1}),\lambda^{j+1}\bigr)\f \bigl(\lambda,\lambda\bigr)\f^{\frac{j(j-1)}{2}}
\bigl(\bigl(x((gh)_{_1})^{i(\lambda^{-1})}\bigr),x((g^{-1})_{_1})\bigr)\f \\ \nonumber
&=&\bigl(x(g_{_1})x(h_{_1}),\lambda^{j+1}\bigr)\bigl(\lambda,\lambda\bigr)\f^{\frac{j(j-1)}{2}}
\bigl(\lambda^j x(g_{_1})x(h_{_1}),x(g_{_1})\bigr)\f.\end{eqnarray}
\end{proof}
\begin{proposition} \label{mgsp para} $$\mgpt= \overline{\gmt} \ltimes \ut,$$
$$\! \! \mgppt= \mgmt \ltimes \ut.$$
\end{proposition}
\begin{proof} Since $\mgpt=i(\F^*)\mmt$ and since $\mmt$ normalizes $\ut$ it is sufficient to prove that for $n \in \ut$ and $\lambda \in \F^*$, $$\bigl(i(\lambda),1\bigr)n\bigl(i(\lambda),1\bigr)^{-1}=\bigl(i(\lambda)n(i(\lambda^{-1}),1 \bigr).$$
This follows directly from \eqref{the fact}.
\end{proof}
\begin{proposition} \label{+} Let $\mgpmt$ be a standard parabolic subgroup of $\mgpspn$. The following hold:\\
1. $Z\bigl(\mgpmt \bigl)=\overline{Z \bigl(\gpmt \bigr)}$. In particular, $Z \bigl(\mgpspn \bigr)=\overline{\F^* I_{2n}}$. \\
2. The relation between $\mgpmt$ and $\mmt$ is similar to the relation between $\gpmt$ and $\mt$. More precisely, \begin{equation} \label{sp gsp+} \mgpmt=Z \bigl(\mgpspn \bigr)  \mmt \end{equation} and $Z\bigl(\mgpspn \bigr)  \cap \mmt=\overline{\pm I_{2n}}.$ \\
3. Two elements in $\mgpmt$  commute if and only if their projections to $\gpmt$ commute.\\
4. Let $\pi$ be a representation of $\mgpmt$. Then $\pi$ is irreducible if and only if its restriction to $\mmt$ is irreducible.\\
5. Let $\tau$ be an irreducible smooth admissible representation of $\mmt$. Let $\chi$ be a character of $Z \bigl(\mgpspn \bigr)$ which agrees with $\chi_\tau$ on $\overline{\pm I_{2n}}$. For $g \in \mgpmt$ write $g=zg_{0}$, where $z \in Z \bigl(\mgpspn \bigr)$ and $g_{_0} \in \mmt$. The map
$$g \mapsto \chi(z)\tau(g_{_0})$$ defines an irreducible representation of $\mgpmt$ and every irreducible smooth admissible representation of $\mgpmt$ is obtained that way.\\

\end{proposition}
\begin{proof} Clearly $Z \bigl(\mgpmt \bigr)  \subseteq \overline{Z \bigl(\gpmt \bigr)}$. Suppose now that $g \in Z \bigl(\gpmt \bigr)$. Then, by \eqref{the fact} $$(g,\epsilon)(h,\epsilon')(g,\epsilon)^{-1}=(h,\epsilon')$$ for any $h \in \gpmt$. This proves 1. Given 1, 2 is clear. 3 follows from 2 and from the first assertion of Lemma \ref{mspn prop}. The last two assertions follow from the second.
\end{proof}
Note: Let $\tau_{_0}$ be a genuine smooth admissible irreducible representation of $\mspn$. Fix $\sigma \in \mgspn$ such that $\lambda(\sigma)=-1$. It was proven by Sun in Theorem 1.2 of \cite{SB} that
\begin{equation} \label{sun dual} \widehat{\tau_{_0}} \simeq \tau_{_0}^\sigma. \end{equation}
It now follows from \eqref{+} above that if $-1 \in {\F^*}^2$ then any genuine smooth admissible irreducible representation of $\mspn$ is self dual.

{\bf  Definition}: Let $\bold{t} = (n_1,n_2,\ldots,n_r;n_{r+1})$. We say that $\mt$, $\gpmt$, $\gmt$, $\pt$, $\gppt$, $\gpt$ or their metaplectic double covers are of odd type if at least one of the numbers $n_1,n_2,\ldots,n_r,n_{r+1}$ is odd. Otherwise we say that these groups are of even type. Thus, $\gspn$, $\gsppn$ and $\spn$ are of odd type if and only if $n$ is odd. Furthermore, provided that $n$ is odd, all the standard parabolic subgroups of $\gspn$, $\gsppn$ and $\spn$ are of odd type. Note that, regardless of the parity of $n$, the Borel subgroups, $\bspn$, $\bgsppn$ and $\bgspn$ are of odd type. This definition is motivated by the following crucial observation:
\begin{lem} \label{x on z} Let $\gpmt$ be a standard parabolic subgroup of $\gsppn$. Define
$$Z_{_1} \! \bigl(\gpmt \bigr)=\bigl \{g_{_1} \! \! \mid g \in Z\bigl(\gpmt \bigr) \bigr \}.$$
If $\gpmt$ is of odd type then
$$x\bigl(Z_{_1} \! \bigl(\gpmt \bigr)\bigr)={\F^*}/{\F^*}^2$$
while if $\gmt$ is of even type then
$$x\bigl(Z_{_1} \! \bigl(\gpmt \bigr)\bigr)=1.$$
\end{lem}
\begin{proof} Assume that $\bold{t}=(n_1,n_2,\ldots,n_r;n_{r+1})$. An element $g\in Z \! \bigl(\gpmt \bigr)$ has the form
$$g=diag \bigl(a_1I_{n_1},a_2I_{n_2},\ldots,a_rI_{n_r},bI_{n_{r+1}},
b^2a_1^{-1}I_{n_1},b^2a_2^{-1}I_{n_2},\ldots,b^2a_r^{-1}I_{n_r},bI_{n_{r+1}} \bigr),$$
where $a_1,a_2,\ldots,a_r,b \in \F^*$. Since $\lambda(g)=b^2$ we have
$$g_{_1}=diag \bigl(a_1I_{n_1},a_2I_{n_2},\ldots,a_rI_{n_r},bI_{n_{r+1}},
a_1^{-1}I_{n_1},a_2^{-1}I_{n_2},\ldots,a_r^{-1}I_{n_r},b^{-1}I_{n_{r+1}}\bigr).$$
Thus, by \eqref{x def}, $$x(g_{_1})=  b^{n_{r+1}}\prod_{n_k \in \N_{odd}}    a_k.$$
\end{proof}

\begin{proposition} \label{even odd basic} Let $\mgpmt$ be a standard Levi subgroup of $\gspn$. The following hold:\\
1. For $h \in Z \bigl(\mgpmt \bigr) $ and $g \in \mgmt$ we have

\begin{equation} \label{peculiar} (h,\epsilon)(g,\epsilon')(h,\epsilon)^{-1}=\bigl(g,d(g,h)\epsilon'\bigr),\end{equation} and
\begin{equation} \label{peculiar2} (g,\epsilon)(h,\epsilon')(g,\epsilon)^{-1}=\bigl(h,d(g,h)\epsilon'\bigr),\end{equation}
where
$$d(g,h)=\bigl(x(h_{_1}),\lambda(g)\bigr)\f.$$
2. $Z \bigl( \mgmt \bigr)$ consists of elements of the form $(h,\epsilon)$ where $h \in Z\bigl(\gpmt \bigr)$ is such that $$x\bigl(h_{_1}\bigr)\in {\F^*}^2.$$
In particular $Z \bigl( \mgmt \bigr) \subseteq Z \bigl( \mgpmt \bigr)$.\\
3. If $\mgmt$ is of even type then $Z \bigl( \mgmt \bigr) = Z \bigl( \mgpmt \bigr)$.\\
4. If $\mgmt$ is of odd type then $Z \bigl(\mgpmt \bigr) /Z \bigl( \mgmt \bigr) \simeq \F^*/{\F^*}^2$. An isomorphism is given by $(g,\epsilon) \mapsto x(g_{_1})$.\\
5. For $(g,\epsilon),(g', \epsilon') \in Z \bigl(\mgpmt \bigr) $ we have \begin{equation} \label{center structure}(g,\epsilon)(g', \epsilon')=\bigl(gg',\epsilon\epsilon'(x(g_{_1}),x(g'_{_1}))\f\bigr).\end{equation}
 In particular, $\mgspn$ splits over $Z \bigl( \mgmt \bigr)$ via the trivial section.
\end{proposition}
\begin{proof}
\eqref{peculiar}  follows from \eqref{the fact}. To prove \eqref{peculiar2} note that by \eqref{extension}
$$\wc \bigl(i(y),s \bigr)=1,$$
for all $y \in \F^*,$ and $s \in \gspn$. Thus,
\begin{eqnarray} \nonumber (g,\epsilon)(h,\epsilon')(g,\epsilon)^{-1} &=&
(i(\lambda(g)),1)(g_{_1},\epsilon)(h,\epsilon')(g_{_1},\epsilon)^{-1}(i(\lambda(g)),1)^{-1} \\ \nonumber
&=&(i(\lambda(g)),1)(h,\epsilon')(i(\lambda(g)),1)^{-1}.\end{eqnarray}
\eqref{peculiar2} follows now from \eqref{the fact}. We now prove the second assertion of this Lemma. Fix $h\in \gmt$ such that $\lambda(h)$ is not a square. We may choose $g \in T_{2n}(\F) \subset \gmt$ such that $$\bigl(x(g),\lambda(h) \bigr)\f=-1.$$
By \eqref{the fact} $$(g,1)(h,\epsilon)(g,1)^{-1}=(ghg^{-1},-\epsilon).$$
Thus, $(h,\epsilon)$ and $(g,1)$ do not commute. This implies that $Z \bigl( \mgmt \bigr) \subseteq Z \bigl(\mgpmt \bigr)$. The second assertion follows now from the first. The third and forth assertions follow from the second assertion and from Lemma \ref{x on z}. The last assertion follows from the cocycle formula, \eqref{extension}.
\end{proof}
From this point we denote $Z \bigl(\mgpmt \bigr) /Z \bigl( \mgmt \bigr)$ by $\zt$. If $\mgmt \subseteq \mgspn$, where $n$ is odd, then as a set of representatives of $\zt$ we may take $$\{\bigl(bI_{2n},1\bigr)\mid b \in \F^* /{\F^*}^2 \}.$$
The following is a peculiar result of Proposition \ref{even odd basic}:
\begin{proposition} \label{strange} Let $\mgmt$ be a parabolic subgroup of $\mgspn$. Assume that  $\mgmt$ is of odd type. Let $\pi$ be a genuine representation of $\mgmt$ and let $\chi$ be a quadratic character of $\F^*$. Then,
 $$(\chi \! \circ \! \lambda) \otimes \pi \simeq \pi.$$
\end{proposition}
\begin{proof}
Fix $a \in \F^* / {\F^*}^2$ such that $\eta_a=\chi$. By the third assertion of Proposition \ref{even odd basic}, one can pick $\sigma=(h,\epsilon) \in Z \bigl(\mgpmt \bigr) $ such that $x(h_{_1})=a$. Clearly, $\pi \simeq \pi^a$. By \eqref{peculiar}, $$g^\sigma=\bigl(I_{2n},(a,\lambda(g))_{_\F}\bigr) g$$ for all $g\in \mgmt$. Since $\pi$ is genuine, $$\pi^\sigma=(\chi \! \circ \! \lambda) \otimes \pi.$$
\end{proof}
There is no linear analog to Proposition \ref{strange}. In fact, it fails already in the $GL_2(\F)$ case.

\section{Relative representation theory} \label{rrt}
\subsection{One dimensional genuine characters}

\begin{lem} \label{odd not iso} Let $\mgpmt$ be a parabolic subgroup of $\mgspn$. Assume that $\mgpmt$ is of odd type. Let $g$ be any element of $\mgmt$ which lies out side $\mgpmt$.\\
1. Let $\chi$ be a genuine character of $Z \bigl(\mgpmt \bigr) $. Then,
$$\chi^g(z,\epsilon)=\eta_{_{\lambda(g)}}\bigl(x(z_{_1})\bigr)\chi(z,\epsilon).$$
In particular $\chi^g \neq \chi$.\\
2. Let $\pi$ be a genuine representation of $\mgpmt$ with a central character. Then $\pi \not  \simeq \pi^{g}$.
\end{lem}
\begin{proof} The first assertion follows directly from \eqref{peculiar} and from Lemma \ref{x on z}. Note that the assumptions that $\chi$ is genuine and that $\mgpmt$ is of odd type are crucial here. The second assertion follows now since $\pi$ and $\pi^g$ have different central characters.
\end{proof}
Remark: Clearly, the first assertion of Lemma \ref{odd not iso} is false for parabolic subgroups of even type. In Theorem \ref{counter ex} we give a counter example, showing that the second assertion of Lemma \ref{odd not iso} is also false for $\mgspn$, where $n$ is even.

Let $\mgmt$ be a standard Levi subgroup of $\mgspn$. Let $\chi$ be a character of $Z \bigl( \mgmt \bigr)$. Denote by $\Omega_\chi$ the set of characters of  $Z \bigl(\mgpmt \bigr) $ which extend $\chi$. By the third and forth assertions of Proposition \ref{even odd basic}, $\Omega_\chi$ is singleton if $\mgmt$ is of even type and $ \# \Omega_\chi=[\F^*:{\F^*}^2]$ if $\mgmt$ is of odd type.
\begin{lem} \label{cha des} Let $\mgmt$ be a standard parabolic subgroup of $\mgspn$. Let $\chi$ be a genuine character of $Z \bigl( \mgmt \bigr)$.\\
1. $\chi$ has the form $$(z,\epsilon) \mapsto \epsilon\eta(z),$$
where $\eta$ is a character of the projection of $Z \bigl( \mgmt \bigr)$ into $\gspn$.\\
2. Fix $\eta'$, a character of $Z \bigl(\gpmt \bigr)$ which extends $\eta$, and fix $\psi$, a nontrivial additive character of $\F$. Then,
\begin{equation} \label{the set}\Omega_\chi=\bigl \{(g,\epsilon) \mapsto  \eta_a \! \cdot \!  \gamma_\psi\bigl(x(g_{_1})\bigr)\epsilon \eta'(g) \mid a \in \F^*/{\F^*}^2 \bigr \}.\end{equation}
3. Any genuine character of $Z \bigl(\mgpmt \bigr) $ has the form
$$(g,\epsilon) \mapsto \epsilon \xi(g) \gamma_\psi\bigl(x(g_{_1})\bigr),$$
where $\xi$ is a character of $Z \bigl(\gpmt \bigr)$ and $\psi$ is a nontrivial additive character of $\F$.
\end{lem}

\begin{proof} The first assertion follows from the fifth assertion of Lemma \ref{even odd basic}. If $\mgmt$ is of even type then the other assertions are trivial. Assume now that $\mgmt$ is of odd type: From \eqref{center structure} it follows that the right hand side of \eqref{the set} indeed contains $[\F^*:{\F^*}^2]$ elements of $\Omega_\chi$. The third assertion follows now from the second assertion and from \eqref{gammaprop}.
\end{proof}
\begin{cor} \label {char of g meta} Any genuine character of $Z \bigl(\mgpspn \bigr)$ has the form
$$(aI_{2n},\epsilon) \mapsto \epsilon \eta'(a)\gamma_\psi(a^n)$$ where $\eta'$ is a character of ${\F^*}$ and $\psi$ is a non-trivial additive character of $\F$. If $n$ odd, any genuine character of $Z \bigl(\mgspn \bigr)$ has the form
$$(aI_{2n},\epsilon) \mapsto \epsilon \eta(a)$$ where $\eta$  is a character of ${\F^*}^2$.
\end{cor}

\subsection{Clifford theory} \label{Clifford theory}
We shall make repeated use of the following result from Clifford theory.
\begin{lem} \label{clifford} Let $G$ be a group and let $H$ be a normal subgroup of finite index of $G$. For $\pi$, a representation $H$, we denote $I(\pi)=Ind_{H}^{G} \pi$.\\
1. $ I(\pi)$ is irreducible if and only if $\pi$ is irreducible and $\pi \not \simeq \pi^g$ for any element $g \in G$ which lies out side H.\\
2. Let $\pi_{_1}$ and $\pi_{_2}$ be two irreducible representations of $H$.
$I(\pi_{_1})\simeq I(\pi_{_2})$ if and only if $\pi_{_2}=\pi_{_1}^g$ for some $g\in G$.
\end{lem}

Let $\mgmt$ be a standard Levi subgroup of $\mgspn$. Assume that $\mgmt$ is of odd type. If $\pi$ is a representation of $\mgmt$ with a central character $\chi_\pi$, we define $\Omega_{\pi}=\Omega_{\chi_\pi}.$ For $\omega \in \Omega_{\pi}$ we define $\pi_\omega$ to be the $\omega$ eigen subspace of $\pi$. It is clearly an $\mgpmt$ invariant subspace. Further more, since $\zt$ is a finite commutative group,
\begin{equation} \label{decomp 1}  \pi\mid_{_{\mgpmt}}=\bigoplus_{\omega \in \Omega_\pi} \pi_\omega. \end{equation}
The projection map \begin{equation} \label{proj} \phi_w:\pi \twoheadrightarrow \pi_\omega \end{equation} is defined by \begin{equation} \label{eigen sur} \phi_w(v)=[F^*:{\F^*}^2]^{-1}  \sum _{b \in \zt} \omega^{-1}(b)\pi(b)(v). \end{equation}
If, in addition, $\pi$ is genuine, then by Lemma \ref{odd not iso}, $\mgmt / \mgpmt$ acts simply and transitively on the set of summands in the right hand side of \eqref{decomp 1}. Thus, if we fix $\omega \in \Omega_\pi$ then
\begin{equation} \label{decomp 2}\pi\mid_{_{\mgpmt}} = \! \! \bigoplus_{a \in \F^*/{\F^*}^2} \! \! \pi_{\omega}^{i(a)}. \end{equation}
By the second assertion of Lemma \ref{odd not iso}, the summands in \eqref{decomp 1} and \eqref{decomp 2} are pairwise non-isomorphic.
\begin{lem} \label{pi is ind} Let $\mgmt$ be a parabolic subgroup of $\mgspn$. Assume $\mgmt$ is of odd type. Let $(\pi,V_\pi)$ be a genuine representation of $\mgmt$ with a central character. Fix $w\in \Omega_\pi$ and denote  $$I(\pi_w)=Ind_{\mgpmt}^{\mgmt} \pi_w.$$
Then, $\pi \simeq I(\pi_w).$
\end{lem}
\begin{proof} For $v \in V_\pi$ define $f(v):{\mgmt} \rightarrow  \pi_w$ by $$g \mapsto \phi_w\bigl(\pi(g)v\bigr).$$
The map $v \mapsto f(v)$ is obviously an $\mgmt$ intertwining map from $V_\pi$ to $I(\pi_w)$. Since $\mgmt$ permute the $\pi_w$'s, this map is injective. For $f \in I(\pi_w)$ define $v(f) \in V_\pi$ by
$$v(f)=\sum_{a \in \F^* / {\F^*}^2} \pi\bigl(i(a),1\bigr)f\bigl(i(a^{-1}),1\bigr).$$
The map $f \mapsto v(f)$ is a $\mgmt$ map from $I(\pi_w)$ to  $V_\pi$. It is left to show that $$f\bigl(v(f)\bigr)=f.$$
Indeed, any  for $y\in \F^*/{\F^*}^2$ we have
\begin{eqnarray} \nonumber  && f\bigl(v(f)\bigr)\bigl(i(y),1\bigr) \\ \nonumber && =[\F^*:{\F^*}^2]^{-1} \! \! \sum_{a \in \F^*/{\F^*}^2}  \sum_{b \in \zt} \omega^{-1}(b)\bigl(\pi(b)\pi(i(y),1)\pi(i(a),1)f\bigr)\bigl((i(a^{-1}),1\bigr) \\ \nonumber
&& =[\F^*:{\F^*}^2]^{-1}\sum_{a \in \F^*/{\F^*}^2} \sum_{b \in \zt} \omega^{-1}(b)\bigl(\pi(i(a),1)\pi(b)^{i(a)}f\bigr)\bigl((i(ya^{-1}),1\bigr)\\ \nonumber
&& =[\F^*:{\F^*}^2]^{-1}\sum_{a \in \F^*/{\F^*}^2} \sum_{b \in \zt} \eta_a(b)\bigl(\pi(i(a),1)f\bigr)\bigl((i(ya^{-1}),1\bigr).\end{eqnarray}
Since $$[\F^*:{\F^*}^2]^{-1} \sum_{b \in \zt} \eta_a(b)=\begin{cases} 1 & a\in {\F^*}^2 \\ 0 & a \not \in {\F^*}^2
\end{cases},$$ we have shown that
$$f\bigl(v(f)\bigr)\bigl(i(y),1\bigr)=f\bigl(i(y),1\bigr).$$
Taking into account the fact that any $f \in I(\pi_\omega)$ is determined on elements of the form $\bigl(i(y),1\bigr)$, where $y\in \F^*/{\F^*}^2$, we are done.
\end{proof}


\begin{proposition} Let $\mgmt$ be a parabolic subgroup of $\mgspn$. Assume that $\mgmt$ is of odd type. Let $(\pi,V_\pi)$ be a genuine representation of $\mgmt$ with a central character. Fix $w\in \Omega_\pi$. Then, $\pi$ is an irreducible $\mgmt$ module if and only if $\pi_w$ is an irreducible $\mmt$ module.
\end{proposition}
\begin{proof} By Lemma \ref{pi is ind} and by Clifford theory, i.e., Lemma \ref{clifford}, we know that $\pi$ is an irreducible $\mgmt$ module if and only if $\pi_w$ is an irreducible $\mgpmt$ module. By the forth assertion of Lemma \ref{+}, $\pi_w$ is an irreducible $\mgpmt$ module if and only if it is an irreducible $\mmt$ module.
\end{proof}
Assume that $\mgmt$ is of odd type. We now define an equivalence relation on the set of isomorphism classes of smooth admissible  irreducible genuine representations of $\mgpmt$: We say that two smooth admissible irreducible genuine representations of $\mgpmt$, $\pi$ and $\pi'$ are $\mgmt$ equivalent if $\pi^g \simeq \pi'$ for some $g\in \mgmt$. Note that for any $\pi$, an irreducible smooth admissible genuine representation of $\mgpmt$, the $\mgmt$ equivalence class of $\pi$ is of cardinality $[\F^*:{\F^*}^2]$.  We denote the factor set by $E_{\bold{t}}(\F)$. We have proven the following:

\begin{thm} \label{corr rep}There is a one to one onto map from $E_{\bold{t}}(\F)$ to the set of smooth admissible genuine irreducible representations of $\mgmt$. It is given by $$\pi \mapsto Ind_{\mgpmt}^{\mgmt}\pi.$$ Its inverse is
$$\pi \mapsto \pi_\omega,$$ where $\omega$ is  any element of $\omega_\pi$.
\end{thm}

For $\mgmt=\overline{GSp_2(\F)}=\overline{GL_2(\F)}$, these results were proven by Gelbart and Piatetski-Shapiro in \cite{GP}. If $\mgmt=T\overline{'_{2n}(\F)}$, then Theorem \ref{corr rep}
coincides with Theorem 5.1 and Corollary 5.2 of \cite{Mc}. The crucial facts in \cite{Mc}, which deals with central extensions of split reductive algebraic groups over a non-archimedean
local fields, is that $\overline{T'_{2n}(\F)}$ is a Heiznberg group, i.e, a two step nilpotent group, whose center is of finite index and that $\overline{T^{^+}_{2n}(\F)}$ is a maximal commutative subgroup of $\overline{T'_{2n}(\F)}$.

For $n$ even, it is not true that $\pi \not \simeq \pi^g$ for all irreducible $\mgpspn$ modules $\pi$ and $g \in \mgspn$ which lies outside $\mgpspn$. For odd $n$, there is no multiplicity one when restricting a genuine irreducible representation of $\mgspn$ to $\mspn$. Both counter examples is given by certain unitary principle series representations, see Theorem \ref{counter ex}.

For the n-fold covering groups of $\gln$ and $SL_n(\F)$ an analog of Theorem \ref{corr rep} was proven by Adams in \cite{A} ,using the same argument. In fact, an analog of Theorem 3.3 of \cite{A} holds in our context as well. We cite it here but omit its proof since it goes word for word as Adam`s proof.
\begin{thm}\label{adams} Let $\mgpmt$ be a Levi subgroup of odd type. Let $\pi$ be a smooth admissible genuine representation of
$\mgpmt$ with a central character. Denote $$\Pi=Ind_\mgpmt^\mgmt \pi.$$
Assume the $\pi$ has a character, then $\Pi$ has also a character. If we denote these characters by $\Theta_\pi$ and $\Theta_\Pi$ respectively then \begin{equation}\label{adams formula}
\Theta_\pi(g)=[F^*:{\F^*}^2]^{-1}\sum _{b \in \zt} {\chi_{_\pi}}^{\! \! \! \! \! -1}(b)\Theta_\Pi(bg). \end{equation}
\eqref{adams formula} remains true if $\Theta_\pi$ is replaced by $\Theta_{\pi_{_0}}$, where $\pi_{_0}$ is the restriction of $\pi$ to $\mmt$.
\end{thm}

As an immediate application of Theorem \ref{corr rep} we use \eqref{sun dual} to prove the following.


\begin{thm} \label{odd dual} Assume that $n$ is odd. Let $\pi$ be a genuine smooth admissible irreducible representation of $\mgspn$. Let $\eta'$ be a character of $\F^*$ such that
$$\omega_\pi(aI_{2n},\epsilon) =\epsilon\eta(a)$$
for all $a \in {\F^*}^2$ ($\eta'$ is determined up to quadratic twists, see Corollary \ref{char of g meta}). We have:
\begin{equation} \label{dual iso}\widehat{\pi} \simeq  ({\eta'}^{-1} \! \circ \! \lambda) \otimes \pi.\end{equation}
In particular, any genuine smooth admissible irreducible representation of $\mgspn$ whose central character is $(aI_{2n},\epsilon) \mapsto \epsilon$ is self dual. Furthermore, any  genuine smooth admissible irreducible representation of $\mgspn$ is a one dimensional non-genuine twist of a self dual genuine smooth admissible irreducible representation.
\end{thm}
Note that \eqref{dual iso} resembles the well known formula for $GL_2(\F)$.
\begin{proof} By Theorem \ref{corr rep} and Corollary \ref{char of g meta}, we have $$\pi=Ind_{\mgpspn}^{\mgspn} \tau,$$
where $\tau$ is a genuine, smooth admissible irreducible representation of $\mgpspn$ whose central character  has the form $$(aI_{2n},\epsilon) \mapsto \eta'(a)\gamma_\psi(a).$$ Fix $\sigma \in \mgspn$ such $\lambda(\sigma)=-1$. In Lemma \ref{dual+} below we prove that
$$\widehat{\tau} \simeq  ({\eta'}^{-1} \! \circ \! \lambda) \otimes \tau^\sigma.$$
Thus, we have

$$\widehat{\pi} \simeq Ind_{\mgpspn}^{\mgspn} \widehat{\tau} \simeq  Ind_{\mgpspn}^{\mgspn} ({\eta'}^{-1} \! \circ \! \lambda) \otimes \tau^\sigma.$$
Note that if $\chi$ is a character of $\F^*$ and $\varpi$ is a $\mgpmt$ module then
\begin{equation} \label{sim twist}\phi:Ind_{\mgpmt}^{\mgmt}(\chi \! \circ \! \lambda)\otimes \varpi \rightarrow (\chi \! \circ \! \lambda) \otimes Ind_{\mgpmt}^{\mgmt} \varpi\end{equation} defined by

$$(\phi f)(h)=(\chi \! \circ \! \lambda)(h)f(h)$$
is an $\mgmt$ isomorphism. This implies that $$\widehat{\pi} \simeq ({\eta'}^{-1} \! \circ \! \lambda) \otimes Ind_{\mgpspn}^{\mgspn}  \tau^\sigma.$$
Using Theorem \ref{corr rep} again, the first assertion of this theorem follows. The other assertions are clear.
\end{proof}
\begin{lem} \label{dual+} Let $\tau$ be a genuine smooth admissible irreducible representation of $\mgpspn$. By Corollary \ref{char of g meta}, the central character of  $\pi$ has the form $$(aI_{2n},\epsilon) \mapsto \epsilon \eta'(a)\gamma_\psi(a^n)$$
where $\chi$ is a character of  $\, \F^*$ and $\psi$ is a non-trivial additive character of $\F$. Fix $\sigma \in \mgspn$ such $\lambda(\sigma)=-1$. Then, $$\widehat{\tau} \simeq  ({\eta'}^{-1} \! \circ \! \lambda) \otimes \tau^\sigma$$
\end{lem}
\begin{proof} By \eqref{sun dual} and \eqref{sp gsp+} it is enough to show that
$$({\eta'}^{-1} \! \circ \! \lambda) \otimes \tau^\sigma(aI_{2n},\epsilon)=\omega_\tau^{-1}(aI_{2n},\epsilon)$$
for all $a \in \F^*$. If $n \in \N_{odd}$ then by the first assertion of Lemma \ref{odd not iso} $$({\eta'}^{-1} \! \circ \! \lambda) \otimes \tau^\sigma(aI_{2n},\epsilon)=\epsilon \eta'^{-1}(a^2) \eta'(a) \gamma_\psi(a)(a,-1)\f.$$
By \eqref{gammaprop} we are done in this case. The even case follows immediately, since in this case $$\tau^\sigma(aI_{2n},\epsilon)=\tau(aI_{2n},\epsilon)=\epsilon \eta'(a).$$
\end{proof}

\subsection{Parabolic induction}
All parabolic inductions are assumed to be normalized and, as usual, the inducing representation is assumed to be trivial on the unipotent radical. Under this assumption, we identify the representations of a parabolic subgroup with the representations of its Levi part. Notation: Suppose $\gmto \subseteq \gmtt$ are two standard Levi subgroups of $\gspn$. Define
\begin{eqnarray} \nonumber \gptot &=&\gpto \cap \gmtt,\\ \nonumber
\gpptot &=&\gppto \cap \gpmtt=\gsppn \cap \gptot,\\ \nonumber
\ptot &=&\pto \cap \mtt= \spn \cap \gptot.\end{eqnarray}

\begin{thm} \label{para thm} Let $\mgmto \subseteq \mgmtt$ be two Levi subgroups of $\mgspn$. Assume that $\mgmto$ is of odd type. Let $\pi$ be a genuine smooth admissible irreducible representation of $\mgmto$. Let $\tau$ be any irreducible $\mmto$ module which appears in the restriction of $\pi$ to $\mmto$.  Define
$$\Pi=Ind_{\mgptot}^{\mgmtt}\pi$$
and
$$\Delta=Ind_{\mptot}^{\mmtt}\tau.$$
Then,\\ \\
1. If $\mgmtt$ is of odd type then $\Pi$ is irreducible if and only if $\Delta$ is irreducible.\\
2. If $\mgmtt$ is of even type and $\pi$ is supercuspidal then $\Pi$ is irreducible if and only if $\Delta$ is irreducible
and $\tau$ is not Weyl conjugate to $\tau^{g}$ for any $g \in \mgmto$ which lies outside $\mgpmto$.
\end{thm}

\begin{proof}  Let $\tau'$ be any irreducible $\mgpmto$ module which appears in the restriction of $\pi$ to $\mgpmto$ such that its restriction to $\overline{M_{\overrightarrow{t_1}}}$ is (isomorphic to) $\tau$. By Theorem \ref{corr rep} $$\pi=Ind_{\mgpmto}^{\mgmto} \tau'.$$
Thus, by using induction by stages twice we get
\begin{equation} \label{ind stg}\Pi \simeq Ind_{\mgpptot}^{\mgmtt}\tau' \simeq
 Ind_{\mgpmtt}^{\mgmtt}\Delta',\end{equation}
where $$\Delta'=Ind_{\mgpptot}^{\mgpmtt}\tau'.$$

Note that since $Z \bigl(\mgpmtt \bigr) \subset Z \bigl(\mgpmto \bigr)$ and since $\mgpmtt=Z \bigl(\mgpmtt \bigr)\mmtt$, the map
$$f\mapsto f \! \mid_{_{\mmtt}}$$ is an $\mmtt$ isomorphism from
$\Delta' \! \mid_{_{\mmtt}}$ to $\Delta$.

Assume first that $\mgmtt$ is of odd type. By Theorem \ref{corr rep}, $\Pi$ is irreducible if and only if $\Delta'$ is irreducible.
On the other hand, by Proposition \ref{+}, $\Delta'$ is irreducible if and only if its restriction to  $\mmtt$ is irreducible.

Assume now that $\mgmtt$ is of even type. By Lemma \ref{clifford}, i.e, clifford theory, $\Pi$ is irreducible if and only if $\Delta'$ is irreducible and $\Delta' \not \simeq {\Delta \!'}^{^{i(a)}}$ for any $a \in \F^*$ which is not a square. As explained above, $\Delta'$ is a reducible $\mgpmtt$ module if and only if $\Delta$ is a reducible $\mmtt$ module. Furthermore, since $\mgmtt$ is of even type, it follows from Proposition \ref{even odd basic} that $\Delta'$ and ${\Delta \! '}^{^{i(a)}}$ have the same central character. Thus,  $\Delta'  \simeq {\Delta \! '}^{^{i(a)}}$ if and only if $\Delta \simeq {\Delta}^{i(a)}$. Note that \begin{equation} \label{twist iso}  \Delta^{i(a)} \simeq Ind_{\mptot}^{\mmtt}\tau^{i(a)}.\end{equation}
Indeed, an isomorphism is given by  $f \mapsto \widehat{f}$, where $\widehat{f}(g)=
f(g^{i(a)})$. Thus, if $\tau \simeq \tau^{i(a)}$ then  $\Delta \simeq {\Delta}^{i(a)}$. Since  $\tau$ is supercuspidal, it follows from Theorem 2.9 of \cite{BZ77} which extends in a straight forward way to the metaplectic case, that given that $\Delta$ is irreducible then $\Delta \simeq {\Delta}^{i(a)}$ if and only if $\tau$ is Weyl conjugate to $\tau^{i(a)}$.
\end{proof}
We now show that in the even cases, analogs to Theorem \ref{dual+} and Proposition \ref{strange} hold in the context of a parabolic induction, provided that the inducing parabolic group is of odd type.
\begin{proposition} \label{even dual} Let $\mgpt$ be a parabolic subgroup of $\mgspn$. Assume that $\mgpt$ is of odd type and that $n$ is even. Let $\pi$ be a genuine smooth admissible irreducible representation of $\mgmt$. Define
$$\Pi=Ind_{\mgpt}^{\mgspn}\pi.$$ If $\Pi$ is irreducible then
$$\widehat{\Pi} \simeq  ({\eta}^{-1} \! \circ \! \lambda) \otimes \Pi,$$
where $\eta$ is the character of $\F^*$ defined by the relation
$$\chi_{_\Pi}(aI_{2n},\epsilon)=\chi_{\pi}(aI_{2n},\epsilon)=\epsilon \eta(a)$$
for all $a \in {\F^*}$, see Corollary \ref{char of g meta}.

\end{proposition}
\begin{proof} Fix $\omega \in \Omega_\pi$. We have  $$\Pi \simeq Ind_{\mgpspn}^{\mgspn}\Delta',$$
where $$\Delta'=Ind_{\mgppt}^{\mgpspn} \pi_\omega.$$
Since $n$ is even $Z \bigl(\mgspn \bigr)=Z \bigl(\mgpspn \bigr).$ This implies that $$\chi_{_{\Delta'}}=\chi_{_\Pi}=\chi_\pi \! \mid_{_{Z (\mgspn)}}.$$
The irreducibility of $\Pi$ implies that $\Delta'$ is irreducible. Thus, by Lemma \ref{dual+} we have
\begin{eqnarray} \nonumber \widehat{\Pi} &\simeq& Ind_{\mgpspn}^{\mgspn}\widehat{\Delta'} \simeq
Ind_{\mgpspn}^{\mgspn} ({\eta}^{-1} \! \circ \! \lambda) \otimes {\Delta \!'}^{^{i(-1)}} \\ \nonumber &\simeq&
({\eta}^{-1} \! \circ \! \lambda) \otimes  Ind_{\mgpspn}^{\mgspn} {\Delta \!'}^{^{i(-1)}} \simeq
({\eta}^{-1} \! \circ \! \lambda) \otimes \Pi.\end{eqnarray}
\end{proof}
Note that we can replace the condition the $\Pi$ is irreducible with the condition that $\Delta'$ is irreducible. Unlike to odd case, this is a weaker condition.
\begin{proposition}  Let $\mgmto \subseteq \mgmtt$ be two Levi subgroups of $\mgspn$. Assume that $\mgmto$ is of odd type. Let $\pi$ be a genuine representation of $\mgmto$. Define
$$\Pi=Ind_{\mgptot}^{\mgmtt}\pi.$$
Let $\chi$ be a quadratic character of $\F^*$. Then,
 $$(\chi \! \circ \! \lambda) \otimes \Pi \simeq \Pi.$$
\end{proposition}
\begin{proof} By proposition \ref{strange} and by \eqref{sim twist},
$$\Pi=Ind_{\mgptot}^{\mgmtt}\pi \simeq Ind_{\mgptot}^{\mgmtt} (\chi \! \circ \! \lambda) \otimes \pi \simeq (\chi \! \circ \! \lambda) \otimes Ind_{\mgptot}^{\mgmtt}  \pi=(\chi \! \circ \! \lambda) \otimes \Pi.$$
\end{proof}
\section{Irreducibility theorems} \label{irr thm}
\begin{thm} Assume that $n$ is odd. Let $\pi$ be a unitary smooth admissible irreducible genuine supercuspidal representation of $\overline{P'_{n;0}(\F)}$. Then, $$Ind_{\overline{P'_{n;0}(\F)}}^{\mgspn} \pi$$ is irreducible.
\end{thm}
\begin{proof} Any $\overline{M_{n;0}(\F)}$ irreducible module which appears in the restriction of $\pi$ to $\overline{M_{n;0}(\F)}$ is clearly a unitary smooth admissible irreducible genuine supercuspidal representation. By Corollary 6.3 of \cite{Sz12}, if $\tau$ is an irreducible genuine supercuspidal representation of  $\overline{M_{n;0}(\F)}$, where $n$ is odd then  $$Ind_{\overline{P_{n;0}(\F)}}^{\mspn} \tau$$ is irreducible. The assertion follows now form Theorem \ref{para thm}.
\end{proof}

\begin{thm} \label{odd unitary} Assume that $n$ is odd. All genuine principal series representations of $\gspn$ which are induced form irreducible genuine unitary representations of $\overline{T'_{2n}(\F)}$ are irreducible.
\end{thm}
\begin{proof} Regardless of the parity of $n>0$, every irreducible genuine representation of $\overline{T'_{2n}(\F)}$ is an $[\F^*:{\F^*}^2]$ dimensional representation induced from a genuine character of $\overline{T^{^+}_{2n}(\F)}$. Thus, if $\pi$ is an irreducible genuine unitary representation of $\overline{T'_{2n}(\F)}$, and $\tau$ is any $\overline{T_{2n}(\F)}$ irreducible module which appears in the the restriction of $\pi$ to $\overline{T_{2n}(\F)}$, then $\tau$ is a unitary character. By Theorem 5.1 of \cite{Sz12}, regardless of the parity of $n$, all the principal series representation of $\mspn$ induced from unitary characters are irreducible. Assuming that $n$ is odd, the assertion follows now from Theorem \ref{para thm}.
\end{proof}
Note that the fact that $n$ is odd is crucial here, see Theorem \ref{counter ex}, for a counter example.

Remark: Let $\mmto \subseteq \mmtt$ be Levi subgroups of $\mspn$. Let $\tau$ be an irreducible genuine representation of $\mmto$. Define $$I(\tau)=Ind_{\mptot}^{\mmtt}\tau.$$ Since $$I\bigl(\tau^{i(a)}\bigr) \simeq I^{i(a)}\bigl(\tau\bigr),$$
It follows that $I\bigl(\tau \bigr)$ is irreducible if and only if $I\bigl(\tau^{i(a)}\bigr)$. In particular, if
$$\mto \simeq GL_{n_1}(\F) \times GL_{n_2}(\F)
\ldots \times GL_{n_r}(\F),$$ then, as explained in Lemma \ref{mspn prop}, $\tau$ has the form
$$\bigl(diag(g_{_1}g_{_2},\ldots,g_{n_r},g_{_1}g_{_2},\ldots,g_{n_r},
\widetilde{g_{_1}},\widetilde{g_{_2}},\ldots,\widetilde{g_{n_r}})\epsilon \bigr) \mapsto \epsilon \gamma_\psi\bigl(\prod_{i=1}^{n_r}\det(g_i)\bigr)\otimes \Big(\otimes_{i=1}^r
{\sigma_i}(g_i)\Bigr),$$
where $\psi$ is a non-trivial character of $\F$ and for $1 \leq i \leq r$, $\sigma_i$ are smooth irreducible representation of $GL_{n_i}(\F)$. In this case $\tau^{i(a)}$ is isomorphic to
\begin{equation} \label{quad twist}\bigl(diag(g_{_1}g_{_2},\ldots,g_{n_r},g_{_1}g_{_2},\ldots,g_{n_r},
\widetilde{g_{_1}},\widetilde{g_{_2}},\ldots,\widetilde{g_{n_r}}),\epsilon \bigr) \mapsto \epsilon \gamma_{\psi_a}\bigl(\prod_{i=1}^{n_r}\det(g_i)\bigr)\otimes \Big(\otimes_{i=1}^r
{\sigma_i}(g_i)\Bigr).\end{equation}
Thus, since $\gamma_{\psi_a}=\eta_a \cdot \gamma_\psi$, the question of reusability of $I(\tau)$ is not sensitive to non-genuine quadratic twists of the inducing representation. This fact has no linear analog: For example, it is well known that $Ind_{B_2(\F)}^{SL_2(\F)} \chi$ is irreducible if $\chi$ is a trivial character but reducible if $\chi$ is a non-trivial quadratic character.

In the case where $\F$ is a p-adic field of odd residual characteristic, using  the results  of Zorn given in \cite{Zo09}, we shall obtain a complete list of reducible principal series representation of $\overline{GSp_{_4}(\F)}$: The principal series representation of $\overline{GSp_{_4}(\F)}$ has the form
$$\pi(\chi)=Ind^{\overline{GSp_{_4}(\F)}}_{\overline{B_4^{^+}(\F)}}\chi,$$
where $\chi$ is a character of $\overline{T_4^{^+}(\F)}$ extended trivially on $N_4(\F)$ to $\overline{B_4^{^+}(\F)}$. By Theorem \ref{para thm}, the question of reducibility of $\pi$ depends only on the restriction of $\chi$ to $\overline{T_4(\F)}$. Any genuine character of $\overline{T_4(\F)}$ has the form
$$g=\bigl(diag(a,b,a^{-1},b^{-1}),\epsilon\bigr) \mapsto (\chi_{_1} \otimes \chi_{_2} \boxtimes \gamma_\psi)(g)=\epsilon \gamma_\psi(ab) \chi_{_1}(a)\chi_{_2}(a),$$
where $\chi_{_1}$ and $\chi_{_2}$ are character of $\F^*$ and $\psi$ is a non-trivial character of $\F$.

\begin{thm} Let $\F$ be a p-adic filed of odd residual characteristic. Let $\chi$ be a genuine character of $\overline{T_4^{^+}(\F)}$. Denote by $$\chi_{_0}=\chi_{_1} \otimes \chi_{_2} \boxtimes \gamma_\psi$$ the restriction of
$\chi$ to $\overline{T_4(\F)}$. Then, $\pi(\chi)$ is reducible if and only if one of the following hold\\
I. $\chi_{_0}$ is Weyl conjugate to $$\chi_{_1} \eta_a \otimes \chi_{_2}\eta_a  \boxtimes \gamma_\psi,$$
where $a$ is a non-square element in $\F^*$.\\
II. $\chi_{_0}$ is Weyl conjugate to $\xi \ab \cdot \ab^{s+\half}  \otimes \xi \ab \cdot \ab^{s-\half} \boxtimes \gamma_\psi$, where $\xi$ is a unitary character of $\F^*$ and $s \in \R$.\\
III. $\chi_{_0}$ is Weyl conjugate to $\xi \ab \cdot \ab^{s}  \otimes \eta_b \ab \cdot \ab^{\half} \boxtimes \gamma_\psi$, where $\xi$ is a unitary character of $\F^*$ and $s \in \R$, $b \in \F^*$.
\end{thm}
\begin{proof} By \eqref{quad twist}, it is clear that $\chi_{_0}$ is Weyl conjugate to $\chi_{_0}^{i(a)}$ for some non-square $a \in \F^*$ if and only if $I$ holds. By Theorem 1.1 of \cite{Zo09}, $Ind^{\overline{Sp_{_4}(\F)}}_{\overline{B_4(\F)}}\chi_{_0}$ is reducible if and only if $II$ or $III$ hold. The assertion follows now from the second part of Theorem \ref{para thm}.
\end{proof}
It should be noted that the results of \cite{Zo09} are given for the $\C^1$ cover of $Sp_{_4}(\F)$ constructed via Leray cocycle. Rao cocycle and Leray cocycle defer by a co-boundary. We have used the explicit formula for this  co-boundary as it appears in Theorem 5.3 of \cite{R} to translate these results to $\overline{Sp_{_4}(\F)}$.

\begin{thm} \label{counter ex} The following counter examples hold:\\
1. Assume that $n$ is odd. There exists an irreducible smooth admissible genuine representation of $\overline{\gspn}$ whose restriction to $\mspn$ is not multiplicity free.\\
2. Assume that $n$ is even. There exists a genuine principal series representation of $\overline{\gspn}$ induced from a unitary representation of $\overline{T'_{2n}(\F)}$ which is reducible.\\
3. Assume that $n$ is even. There exists $\pi$, a genuine irreducible representation of $\mgpspn$
such that $\pi \simeq \pi^g$ for some $g \in \overline{\gspn}$ which lies outside $\mgpspn$.\\
\end{thm}
\begin{proof} We use essentially the same construction to prove all three assertions. Let $\chi$ be a character of $\F^*$. Let $\psi$ be a non-trivial additive character of $\F$. Define $\chi_\psi$ to be the following character of $\overline{T_{2n}(\F)}$:
\begin{equation} \label{counter}\Bigl(\begin{pmatrix} _{a} & _{0}\\_{0} & _{\widetilde{a}}\end{pmatrix},
\epsilon \Bigr) \mapsto   \chi \! \cdot \! \gamma_{\psi}\bigl(\det(a)\bigr)
\epsilon.\end{equation}
Define
$$I \bigl(\chi_\psi \bigr)=Ind^{\mspn}_{\overline{B_{2n}(\F)}} \chi_\psi.$$
Fix $b \in \F^*$ which is not a square. By \eqref{quad twist} and by the same argument we used for \eqref{twist iso}, \begin{equation} \label{one to use} I^{i(b)} \bigl(\chi_\psi \bigr) \simeq I \bigl((\chi \cdot \eta_b) _\psi \bigr). \end{equation}
Suppose now that $I \bigl(\chi_\psi \bigr)$ is irreducible. Then, Using the standard intertwining operator attached to the Weyl element represented by $\omega_n'=\begin{pmatrix}
& _{ }& _{-\omega_n } \\
& _{\omega_n} & _{ }\end{pmatrix},$ where $$\omega_n=\begin{pmatrix}
_{ } & _{ } & _{ } & _{1} \\
_{ } & _{ } & _{1} & _{ } \\
_{ } & _{\upddots} & _{ } & _{ } \\
_{1} & _{ } & _{ } & _{ }
\end{pmatrix},$$
we observe that   \begin{equation} \label{two to use}I \bigl(\chi_\psi \bigr)\simeq I \bigl((\chi^{-1})_\psi \bigr),
\end{equation} see Section 2 of \cite{Sz12}.

Let $\xi$ be any character of $\F^*$ such that $\xi(-1)=\chi^n(-1)$. Extend $I (\chi_\psi )$ to $\mgpspn$ by defining
\begin{equation} \label{extend} (zI_{2n},
\epsilon ) \mapsto \epsilon  \xi(z)
\gamma_{\psi}^{-1}(z^n). \end{equation}
Denote this $\mgpspn$ module by $I_\xi \bigl(\chi_\psi \bigr)$. Note that
$$I_\xi \bigl(\chi_\psi \bigr) \simeq Ind^{\mgpspn}_{\overline{B^{^+}_{2n}(\F)}} \chi_{\psi,\xi},$$ where
$\chi_{\psi,\xi}$ is the character of $\overline{T^{^+}_{2n}(\F)}$ which extends $\chi_\psi$ by \eqref{extend}. Define now $\Pi_{\xi, \chi,\psi}$ to be the following genuine principal series representation of $\mgspn$
$$\Pi_{\xi, \chi,\psi}=Ind_{\mgpspn}^{\mgspn} I_\xi \bigl(\chi_\psi \bigr) \simeq
Ind_{\overline{B'_{2n}(\F)}}^{\mgspn} \biggl(Ind_{B^{^+}_{2n}(\F)}^{B'_{2n}(\F)} \chi_{\psi,\xi} \biggr).$$
From this point we assume that $\chi$ is a character of $\F^*$ of order four. Since the inducing representation is unitary,  $I\bigl(\chi_\psi \bigr)$ is irreducible, see Theorem 5.1 of \cite{Sz12}. Also, $\chi^2$ is a non-trivial quadratic character. This implies that $\chi \cdot \eta_b=\chi^{-1}$ for some non-square $b \in \F^*$. Thus, by \eqref{one to use} and \eqref{two to use} we have $$I\bigl(\chi_\psi \bigr) \simeq I^{i(b)} \bigl(\chi_\psi \bigr).$$

Assume that $n$ is odd. By Theorem \ref{odd unitary}, $\Pi_{\xi, \chi,\psi}$ is an irreducible genuine representation of $\overline{\gspn}$ and both $I \bigl(\chi_\psi \bigr)$ and
$I^{i(b)} \bigl(\chi_\psi \bigr)$ appears in the restriction of $\Pi_{\xi, \chi,\psi}$ to $\mspn$ as disjoint subspaces. This proves the first assertion of the proposition.

Assume now that $n$ is even. By the same argument used in the proof of the second part of Theorem \ref{para thm}, since $I \bigl(\chi_\psi \bigr) \simeq I^{i(b)} \bigl(\chi_\psi \bigr)$ it follows that $i(b)$  fixes $I_\xi \bigl(\chi_\psi \bigr)$. From Lemma \ref{clifford} we conclude that $\Pi_{\xi, \chi,\psi}$ is reducible. This prove the second and third assertions stated here.
 \end{proof}
\begin{cor} Let $\mgmt$ be a parabolic subgroup of $\mgspn$. Assume $\mgmt$ is of odd type. Let $\pi$ be any irreducible smooth admissible genuine representation of $\mgmt$. Contrary to the situation in the linear case, the restriction of $\pi$ to
$\mmt$ is always a sum of $[\F^*:{\F^*}^2]$ irreducible smooth admissible genuine representations of $\mmt$. Furthermore, if $\mgmt=\gspn$, $\mmt=\mspn$ and $n$ is odd, then, unlike the linear case, multiplicity one does not hold.
\end{cor}
\begin{proof} The first assertion follows from Theorem \ref{corr rep} and from the relation between $\mmt$ and $\mgpmt$ proven in Proposition \ref{+}. The second assertion follow from the first counter example given in Theorem \ref{counter ex}.
\end{proof}

Remark: Let $\F$ be a p-adic field of odd residual characteristic. It can be shown the $\mgspn$ splits over $Gsp_{2n}(\Of)$. This splitting is not unique, in fact there are two splittings and the corresponding embeddings are conjugated if and only if $n$ is odd. Fix an embedding $\kappa$ of $Gsp_{2n}(\Of)$ into $\mgspn$. A genuine principal series of $\mgspn$ contains a $\kappa \bigl( Gsp_{2n}(\Of) \bigr)$ vector is and only if it is induced from a character of $\overline{B'_{2n}(\F)}$ whose restriction to $\overline{T^{^+}_{2n}(\Of)}$ is $(g,\epsilon) \mapsto \epsilon\gamma_\psi\bigl(x(g_{_1})\bigr)$. Thus, since $\F$ has unramified character of order 4, the family of examples given in Theorem \ref{counter ex} includes reducible unramified principles series representations of $\mgspn$ where $n$ is even. Furthermore, there is another family of reducible unitary principle series representations which includes  unramified representations: Let $\chi_1$ and $\chi_2$ be quadratic characters of $\F$ and let $\psi$ be a character of $\F$. Define $\chi_\psi$ to be the following character of $\overline{T_{2n}(\F)}$:
$$\Bigl(\begin{pmatrix} _{a} & _{0}\\_{0} & _{\widetilde{a}}\end{pmatrix},
\epsilon \Bigr) \mapsto   \chi_1 \! \cdot \! \gamma_{\psi}\bigl(\det(a)\bigr)\chi_2(a_2 a_4\ldots a_n)
\epsilon. $$
Here $a=diag(a_1 a_2\ldots a_n).$ Pick a character $\xi$ of $\F^*$ which satisfies  $\xi(-1)=\chi_2^{\frac n 2}(-1)$ and extend $\chi$ to $\overline{T^{+}_{2n}(\F)}$ as in \eqref{extend}. Using similar argument to the one used in Theorem \ref{counter ex} it can proven that $\Pi_{\xi, \chi,\psi}$ is an example to 2 and 3 of Theorem \eqref{extend}.

\section{Whittaker functionals} \label{whi}

\subsection{Uniqueness of $\psi \times \omega$-Whittaker functional}

Let $\mt$ be a standard Levi subgroup of $\spn$. Define
$$\nt=\mt \cap N_{2n}(\F).$$
$\nt$ is the unipotent radical of a Borel subgroups of $\mt$, $\gpmt$ and $\gmt$.

Let $\overline{M}$ be one of the groups $\overline{\mt}$, $\mgmt$ or $\mgpmt$. Let $\pi$ be a representation of  $\overline{M}$ and let $\psi$ be a non-degenerate character $\nt$. By a $\psi$-Whittaker functional on $\pi$ we mean a linear functional $\vartheta$ on $V_\pi$ such that
$$\vartheta\bigl(\pi(u,1)v\bigr)=\psi(u)\vartheta(v),$$ for all $v \in V_\pi, \, u \in \nt$. Denote by $W_{\psi ,\pi}$ the space of $\psi$-Whittaker functionals on $\pi$. We say that $\pi$ is $\psi$-generic if $$\dim (W_{\psi ,\pi})>0.$$
For $t \in \overline{T'_{2n}(\F)} \cap \overline{M}$ let $\psi^{(t)}$ be the character of $U_{\bold{t}}(\F)$ defined by $$u \mapsto \psi(tut^{-1}).$$ Clearly, if $\vartheta$ is a $\psi$-Whittaker functional then the functional $\vartheta^{(t)}$ defined by $$v\mapsto \vartheta\bigl(\pi(t)v)\bigr)$$ is a $\psi^{(t)}-$ Whittaker functional on $\pi$ and that $\vartheta \mapsto \vartheta^{(t)}$ is an isomorphism between $W_{\psi ,\pi}$ and $W_{\psi^{(t)} ,\pi^{(t)}}$.

There is only one orbit of non-degenerate characters of $\nt$ under the action of $\overline{T'_{2n}(\F)}$. Thus, any generic representation of $\mgmt$ is generic with respect to all non-degenerate characters. Fix $\bold{t}=(n_1,n_2,\ldots,n_r;n_{r_1})$. If $n_{r+1}=0$ there is also only one orbit of non-degenerate characters of $\nt$ under the action of either $\overline{T_{2n}(\F)}$ or  $\overline{T^{^+} _{2n}(\F)}$. However, if $n_{r+1}>0$ then there are $[\F^*:{\F^*}^2]$ orbits of non-degenerate characters of $N_{\bold{t}}$ under the actions of either  $\overline{T_{2n}(\F)}$ or  $\overline{T^{^+} _{2n}(\F)}$. If we fix $\psi$, a non-degenerate character of $\nt$ , then
\begin{equation} \label{all orbits} \{\psi^{i(a)} \mid a \in \F^* / {\F^*}^2\} \end{equation}
is a set of representatives of all these orbits.
\begin{lem} \label{uni whi} Let $\pi$ be an irreducible smooth admissible representation of either of the groups $\mmt$ or $\mgpmt$ and let $\psi$ be a non-degenerate character of $\nt$. Then, $$\dim (W_{\psi,\pi}) \leq 1.$$
\end{lem}
\begin{proof} If $\pi$ is not genuine then we are reduced to linear groups, where this assertion is well known. We assume now that $\pi$ is genuine. For $\mmt=\mspn$  the assertion was proven in \cite{Sz}. As explained in Lemma \ref{mspn prop}, any irreducible admissible genuine representation of $\mmt$ is essentially a tensor product of
irreducible admissible representations of $GL_{n_i}(\F)$ and an irreducible admissible genuine representation of $\overline{Sp_{n_{r+1}}(\F)}$. Thus, the uniqueness for $\mmt$ follows from the uniqueness for  $GL_{n_i}(\F)$ and from the uniqueness for $\overline{Sp_{n_{r+1}}(\F)}$. The uniqueness for $\mgpmt$ now follows from the simple relation between the representation theory of $\mgpmt$ and the representation theory of $\mgmt$ proven in Proposition \ref{+}.
\end{proof}

Suppose now that $\overline{M}$ is either $\mgmt$ or $\mgpmt$. Let $(\pi,V_\pi)$ be representation of $\overline{M}$, let $\vartheta$ be a $\psi$-Whittaker functional on $\pi$ and let $\omega$ be a character of $Z \bigl(\mgpmt \bigr) $. If, $$\vartheta\bigl(\pi(z)v\bigr)=\omega(z)\vartheta(v),$$ for all $v \in V_\pi,\, z\in Z \bigl(\mgpmt \bigr) $, we say that $\vartheta$ is a $\psi \times \omega$-Whittaker functional on $\pi$. Denote by $W_{\psi \times \omega,\pi}$ the space of $\psi \times \omega$-Whittaker functionals on $\pi$. We say that $\pi$ is $\psi \times \omega$-generic if $$\dim (W_{\psi \times \omega,\pi})>0.$$
Note that if $\overline{M}=\mgpmt$ and if $\pi$ has a central character then any $\psi$-Whittaker functional on $\pi$ is a $\psi \times \chi_\pi$ Whittaker functional and that if $\omega \neq \chi_\pi$ then $\dim (W_{\psi \times \omega,\pi})=0$. From Lemma \ref{uni whi} we now conclude the following.

\begin{lem} \label{uni whi+} Let $\pi$ be a smooth admissible irreducible representation of $\mgpmt$. Then, $$\dim (W_{\psi \times \omega,\pi}) \leq 1.$$
\end{lem}

We need one more simple Lemma before we state and proof our main result for this section.
\begin{lem} \label{last one} Let $\mgmt$ be a standard Levi subgroup of $\mgspn$. Assume that $\mgmt$ is of odd type. Let $\pi$ be a genuine representation of $\mgmt$. $\vartheta \mapsto \vartheta^{i(a)}$ is an isomorphism between $W_{\psi \times \omega, \pi}$ and $W_{\psi^{i(a)} \times \omega^{i(a)}, \pi}$.
\end{lem}
\begin{proof} This follows from the first assertion of Lemma \ref{odd not iso}.
\end{proof}
 \begin{thm} \label{spec uniq} Let $\mgmt$ be a standard Levi subgroup of $\mgspn$. Assume that $\mgmt$ is of odd type. Let $\pi$ be a smooth admissible generic irreducible representation of $\mgmt$. The following hold:\\
1. Fix $\psi$, a non-degenerate character of $\nt$. For any $\omega \in \Omega_\pi$ $$\dim (W_{\psi \times \omega,\pi}) \leq 1$$ and there exists at least one $\omega \in \Omega_\pi$ such that $$\dim (W_{\psi \times \omega,\pi}) = 1.$$
2. Fix $\omega \in \Omega_\pi$.  For any $\psi$, a non-degenerate character of $\nt$, $$\dim (W_{\psi \times \omega,\pi}) \leq 1$$ and there exists at least one $\overline{T_{2n}(\F)}$ orbit of non-degenerate characters of $\nt$ such that $$\dim (W_{\psi \times \omega,\pi})= 1.$$
\end{thm}
\begin{proof} If $\pi$ is not genuine then the Theorem follows from the uniqueness of $\psi$- Whittaker functional for linear groups. We now assume that $\pi$ is genuine.

1. Let $\omega \neq \omega'$ be two elements in $\Omega_\pi$. Any $\vartheta \in W_{\psi \times \omega,\pi}$ vanishes on $\pi_{\omega'}$. Thus, if $\dim \bigl(W_{\psi \times \omega,\pi}\bigr)>1$ then $\dim \bigl(W_{\psi \times \omega,\pi_\omega}\bigr)>1$. By Lemma \ref{uni whi+} this is not possible. This proves uniqueness. We now prove existence:
Pick any non zero element $\vartheta \in W_{\psi,\pi}$. There exists at least one $\omega \in \Omega_\pi$ such that $\vartheta \! \mid_{_{\pi_w}} \neq 0$. Clearly $\vartheta \circ \phi_\omega$ is a $\psi \times \omega$-Whittaker functional on $\pi$. Here $\phi_\omega$ is the projection operator defined in \eqref{eigen sur}.

2. Uniqueness is proven exactly as in 1. To prove existence, fix $\psi,$ a non-degenerate additive character of $\nt$. By 1, there exists $\widetilde{\omega} \in \Omega_\pi$ such that $\dim (W_{\psi \times \widetilde{\omega}, \pi}) \neq 0$. Since $\zt$ permutes the elements of $\Omega_\pi$, it follows from Lemma \ref{last one} that there exists $a\in \F^* / {\F^*}^2$ such that  $\dim (W_{\psi^{i(a)} \times \omega, \pi}) \neq 0$.
\end{proof}

Note: For $\mmt=\overline{GL_2(\F)}$, part 1 of Theorem \ref{spec uniq} was proven by  Gelbart, Howe and Piatetski-Shapiro in \cite{GHP} without using the uniqueness Of Whittaker functional for $\msl$. The uniqueness part of  Theorem \ref{spec uniq} may proven using the machinery of Shalika, Gelfand-Kazhdan and Bernstein-Zelevinsky, see \cite{Sh}, \cite{GK} and \cite{BZ}.

Let $\pi$ be a smooth admissible irreducible genuine generic representation of $\mgmt$. As we shall see, $\dim (W_{\psi,\pi})$ may be greater then 1. However, as explained above, $\dim (W_{\psi,\pi})$ does not depend on $\psi$. Thus, we may denote this number by $W_\pi$. For $\omega \in \Omega_\pi$ define $k_\omega$ to be the number of $\overline{T_{2n}(\F)}$ orbits of non-degenerate characters of $\nt$ such that $\pi_\omega$ is $\psi-$ generic. For $\psi$, a non-degenerate character of $\nt$ denote by $k_\psi$ the number of elements $\omega \in \Omega_\pi$ such that $\pi_w$ is $\psi$-generic. Clearly, both $k_\psi$ and $k_\omega$ are not greater then $[\F^*:{\F^*}^2]$.

\begin{thm} \label{dim}Fix $\bold{t}=(n_1,n_2,\ldots,n_r;n_{r+1})$. Let $\mgmt$ be a Levi subgroup of $\mgspn$. Assume that $\mgmt$ is of odd type. Let $\pi$ be a genuine smooth admissible  irreducible generic  representation of $\mgmt$. Then, $k_\psi$ is independent of $\psi$ and $k_\omega$ is independent of $\omega$. Furthermore, if $n_{r+1}>0$ then $$k_\psi=k_\omega=W_\pi \leq[\F^*:{\F^*}^2]$$
and if $n_{r+1}=0$ then $k_\omega=1$ and $$k_\psi=W_\pi =[\F^*:{\F^*}^2].$$
\end{thm}
\begin{proof} Fix $\psi$, a non-degenerate character of $\nt$. Since $\nt \subset \mgpmt$,
$$(\vartheta,t) \mapsto \vartheta^{(t)}$$
defines a representation of the finite commutative group $\zt$ on $W_{\psi,\pi}$. It was proven in Theorem \ref{spec uniq} that $W_{\psi,\pi}$ decomposes over this group with multiplicity 1. Thus, $\dim W_{\psi,\pi}$ is equal to the number of elements $\omega \in \Omega_\pi$ such that $ \pi_w$ is $\psi \times \omega$ generic.
This shows that $$W_{\pi}= k_\psi \leq [\F^*:{\F^*}^2].$$
If $n_{r+1}>0$ then the equality of $k_\psi$ and $k_\omega$ follows
from \eqref{all orbits} and Lemma \ref{last one}. If $n_{r+1}=0$ then $k_\omega=1$ since there is only one orbit of non-degenerate characters of $\nt$
under the action of  $\overline{T_{2n}(\F)}$. It is left to show that in this case, $$k_\psi=[\F^*:{\F^*}^2].$$ Indeed, since $\pi$ is generic then for at least one $\omega \in \Omega_\pi$, $\pi_\omega$ is $\psi$-generic. By Lemma \ref{last one} and the fact that there is only one orbit of non-degenerate characters of $\nt$
under the action of  $\overline{T_{2n}(\F)}$  we conclude now that for any $\omega \in \Omega_\pi$, $\pi_w$ is $\psi-$generic.
\end{proof}

\begin{cor}
 Let $\mgmt$ be a Levi subgroup of $\mgspn$. Assume that $\mgmt$ is of odd type. Let $\pi$ be a genuine smooth admissible  irreducible generic  representation of $\mgmt$. Let $\Omega_0$ be a subset of $\Omega_\pi$ of maximal cardinality such that $\pi_\omega \simeq \pi_{\omega'}$ as $\mmt$ modules for all $ \omega, \, \omega' \in \Omega_0$. Them $$\# \Omega_0 \leq \dim (W_{\psi,\pi}).$$
\end{cor}
\begin{proof}  This statement is trivial if $n_{r+1}=0$ since in this case, by Theorem \ref{dim}
$$\# \Omega_\pi \leq  [\F^*:{\F^*}^2]= \dim (W_{\psi,\pi}).$$
We assume now that  $n_{r+1}>0$. Fix $\omega \in \Omega_0$ and fix $\psi$ such that $\pi_\omega$ is $\psi$-generic. Clearly, $\pi_{\omega'}$ is $\psi$-generic for all $\omega' \in \Omega_0$. Thus, $k_\psi \geq  \# \Omega_0.$
By theorem \ref{dim} we are done.
\end{proof}
\subsection{Rodier type Heredity}
\begin{lem} \label{her lem} Let $\mgmto \subseteq \mgmtt$ be two Levi subgroups of $\mgspn$. Let $\pi$, $\tau$ and $\sigma$  be smooth admissible irreducible representations of $\mgmto$, $\mmto$ and $\mgpmto$ respectively. Define
$$\Pi=Ind_{\mgptot}^{\mgmtt}\pi,$$
$$\Sigma=Ind_{\mgpptot}^{\mgpmtt}\sigma,$$
$$\Delta=Ind_{\mptot}^{\mmtt}\tau.$$

Let $\Psi$ be a non-degenerate character of $N_{\bold{t_2}}(\F)$ and let $\psi$ be the  non-degenerate character of $N_{\bold{t_1}}(\F)$ obtained by restricting $\Psi$. Then,

\begin{eqnarray} \label{Her G} \dim(W_\Pi) &=& \dim (W_\pi), \\
 \label{Her +}\dim(W_{\Psi,\Sigma})&=& \dim(W_{\psi,\sigma}) \leq 1,\\
 \label{Her S} \dim(W_{\Psi,\Delta})&=& \dim(W_{\psi,\tau}) \leq 1.\end{eqnarray}
\end{lem}
\begin{proof}
All three equalities here are metaplectic analogs to Rodier Heredity for linear groups, see \cite{Rod}. In \cite{Ban}, a similar result, in the context of an n-fold cover of $\gln$ was proven. That proof, which in fact goes along the same lines as Rodier's proof, is written in sufficient generality to apply for our cases also. The two inequalities follows from Lemma \ref{uni whi}. We note that that \eqref{Her S} was an ingredient in the proof of the Uniqueness of Whittaker model for $\mspn$ in \cite{Sz}.
\end{proof}

From Lemma \ref{her lem} and Theorem \ref{dim} we conclude the following.
\begin{cor} Let $\mgmto \subseteq \mgmtt$ be two Levi subgroups of $\mgspn$. Assume that $\mgmto$ is of odd type. Let $\pi$ be an irreducible admissible representation of $\mgmto$. Define
$$\Pi=Ind_{\mgptot}^{\mgmtt}\pi.$$
Then,
$$W_\Pi \leq [\F^*:{\F^*}^2].$$
Further more, if $\bold{t_1}=(n_1,n_2,\ldots,n_r;0)$ and $\pi$ is genuine then
\begin{equation} \label{mc her} W_\Pi =[\F^*:{\F^*}^2].\end{equation}
\end{cor}
Note: For $\Pi$, a genuine principal series representations of $\mgspn$, where $\F$ is a p-adic field of odd residual characteristic, \eqref{mc her} coincides with Theorem 8.1 of \cite{Mc}.

Let $\mgmto \subseteq \mgmtt$ be two parabolic subgroups of $\mgspn$. Assume that $\mgmto$ is of odd type.
Let $\pi$ be a genuine generic irreducible smooth admissible representation of $\mgmto$. Define
$$\Pi=Ind_{\mgptot}^{\mgmtt}\pi.$$
Any $\omega \in \Omega_ \pi$ defines a canonical isomorphism
$$\Pi \simeq Ind_{\mgpptot}^{\mgmtt}\pi_\omega.$$
For $f \in Ind_{\mgpptot}^{\mgmtt}\pi_w$ define $f^{(\omega)}$ to be its restriction to $\mgpmtt$. We define now $\omega$ special $\psi$-Whittaker functional attached to $\Pi$ to be a $\psi$-Whittaker functional on $Ind_{\mgpptot}^{\mgmtt}\pi_w$ which vanish on the the subspace
$$\{f \in Ind_{\mgpptot}^{\mgmtt}\pi_\omega \mid f^{(\omega)}=0 \}.$$

\begin{thm} \label{rodier gen} With the notations and assumptions of the definition above The following hold:\\
1. Fix $\psi$, a non-degenerate character of $N_{\bold{t_2}}(\F)$. For any $\omega \in \Omega_{_\pi}$, the space of $\omega$ special $\psi$-Whittaker functionals attached to $\Pi$ is at most 1 dimensional and there exists at least one $\omega \in \Omega_\pi$ such that this space is non-trivial.\\
2. Fix $\omega \in \Omega_{_\pi}$.  For any $\psi$, a non-degenerate character of $N_{\bold{t_2}}(\F)$ the space of $\omega$ special $\psi$-Whittaker functionals attached to $\Pi$ is at most 1 dimensional and there exists at least one $\overline{T_{2n}(\F)}$ orbit of non-degenerate characters $\psi$ of $N_{\bold{t_2}}(\F)$ such that this space is non-trivial.
\end{thm}
\begin{proof} The map $f \mapsto f^{(\omega)}$ is an $\mgpmtt$ surjection from $Ind_{\mgpptot}^{\mgmtt}\pi_w$ to
$Ind_{\mgpptot}^{\mgpmtt}\pi_w.$ This given an isomorphism between the space of $\omega$ special $\psi$-Whittaker functionals attached to $\Pi$ and the space $\psi$-Whittaker functionals on $Ind_{\mgpptot}^{\mgmtt}\pi_w$. \eqref{Her +} implies now the uniqueness part of this theorem. Fix $\psi$, a non-degenerate character of $N_{\bold{t_2}}(\F)$. Since $\pi$ is generic, by the argument used in the proof of part 1 of Theorem \ref{spec uniq}, there exists at least one $\omega \in \Omega_\pi$ such that $W_{\psi,\pi_{\omega}} \neq 0$. This proves existence in part one of the theorem. The existence in the second part of the Theorem follows from the first part and from Lemma \ref{last one}.
\end{proof}

\begin{thm} \label{rodier} Keep the notations and assumptions above. Assume in addition that $\mgmtt$ is of odd type. The following hold:\\
1. For any $\omega \in \Omega_{_\Pi}$ there exists a unique element $\widetilde{\omega} \in \Omega_\pi$ which extends $\omega$. Realize $\Pi$ as $Ind_{\mgpptot}^{\mgmtt}\pi_{\widetilde{\omega}}.$
The space of $\widetilde{\omega}$ special $\psi$-Whittaker functionals attached to $\Pi$ is $W_{\psi \times \omega,\Pi}$.\\
2. Fix $\psi$, a non-degenerate character of $N_{\bold{t_2}}(\F)$. For any $\omega \in \Omega_{_\Pi}$ $$\dim (W_{\psi \times \omega,\Pi}) \leq 1$$ and there exists at least one $\omega \in \Omega_\Pi$ such that $$\dim (W_{\psi \times \omega,\Pi}) = 1.$$
3. Fix $\omega \in \Omega_{_\Pi}$.  For any $\psi$, a non-degenerate character of $N_{\bold{t_2}}(\F)$. $$\dim (W_{\psi \times \omega,\Pi}) \leq 1$$ and there exists at least one $\overline{T_{2n}(\F)}$ orbit of non-degenerate characters $\psi$ of $N_{\bold{t_2}}(\F)$ such that $$\dim (W_{\psi \times \omega,\Pi}) = 1.$$
\end{thm}
\begin{proof}. We prove only the first part of the Theorem. The second part and the third part of this Theorem follow from the first part and from Theorem \ref{rodier gen}. From \eqref{the set} it follows that a bijection between  $\Omega_\pi$ and $\Omega_\Pi$ is given by restriction. This explains the existence and uniqueness of $\widetilde{\omega}$ above. It is left to show that $$\Pi_\omega \simeq Ind_{\mgpptot}^{\mgpmtt}\pi_{\widetilde{\omega}}.$$
We compute $\phi_\omega(\Pi)=\Pi_\omega$:

$$(\phi_\omega f)(g)=[F^*:{\F^*}^2]^{-1}  \! \! \sum _{b \in \ztt} \! \! \omega^{-1}(b)f(gb)=
[F^*:{\F^*}^2]^{-1}f(g)  \! \!\sum _{b \in \ztt}  \! \! \eta_{_{\lambda(g)}}\bigl(x(b_{_1})\bigr).$$
Thus, $\Pi_\omega$ consists of the functions which are supported on $\mgpmto$.
\end{proof}

\end{document}